\documentclass[10pt]{amsart}
\usepackage{pb-diagram,times,amssymb,epsfig,tikz}
\usetikzlibrary{petri}

\def\O{\mathcal O}
\def\Q{\mathbb Q}
\def\C{\mathbb C}
\def\H{\mathfrak H}
\def\P{\mathbb P}
\def\HH{\mathbb H}
\def\R{\mathbb R}
\def\Z{\mathbb Z}

\def\SL{\mathrm{SL}}

\def\tr{\mathit{tr}}
\def\nm{\mathit{n}}
\def\JS#1#2{\left(\frac{#1}{#2}\right)}
\def\M#1#2#3#4{\begin{pmatrix}#1&#2\\#3&#4\end{pmatrix}}
\def\SM#1#2#3#4{\left(\begin{smallmatrix}#1&#2\\#3&#4\end{smallmatrix}\right)}
\def\gen#1{\langle #1\rangle}
\def\gauss#1{\left\lfloor #1\right\rfloor}

\newcounter{lemma}
\newtheorem{Theorem}{Theorem}
\newtheorem{Lemma}[lemma]{Lemma}
\newtheorem{Proposition}[lemma]{Proposition}

\theoremstyle{definition}
\newtheorem{Remark}[lemma]{Remark}
\newtheorem{Definition}[lemma]{Definition}

\begin{document}
\title[Hypergeometric functions]{Algebraic transformations of
  hypergeometric functions and automorphic forms on Shimura curves}

\author{Fang-Ting Tu}
\address{Department of Applied Mathematics, National Chiao Tung
University, 1001 Ta Hsueh Road, Hsinchu, Taiwan 30010, ROC}
\email{ft12.am95g@nctu.edu.tw}
\author{Yifan Yang}
\address{Department of Applied Mathematics, National Chiao Tung
University and National Center for Theoretical Sciences, Hsinchu 300, Taiwan}
\email{yfyang@math.nctu.edu.tw}
\date\today
\thanks{The authors were partially supported by Grant
  99-2115-M-009-011-MY3 of the National Science Council, Taiwan (R.O.C.).}

\begin{abstract} In this paper, we will obtain new algebraic
  transformations of the ${}_2F_1$-hypergeometric functions. The main
  novelty in our approach is the interpretation of identities among
  ${}_2F_1$-hypergeometric functions as identities among automorphic
  forms on different Shimura curves.
\end{abstract}

\maketitle

\begin{section}{Introduction}
For a real number $a$ and a nonnegative integer $n$, let
$$
  (a)_n=\begin{cases}
  1, &\text{if }n=0, \\
  a(a+1)\ldots(a+n-1), &\text{if }n\ge 1, \end{cases}
$$
be the Pochhammer symbol. Recall that, for real numbers $a,b,c$ with $c\neq 0,-1,-2,\ldots$, the ${}_2F_1$-hypergeometric function is defined by
$$
  {}_2F_1(a,b;c;z)=\sum_{n=0}^\infty\frac{(a)_n(b)_n}{(c)_nn!}z^n
$$
for $z\in\C$ with $|z|<1$. The hypergeometric function is a solution
of the differential equation
$$
  \theta(\theta+c-1)F-z(\theta+a)(\theta+b)F=0, \qquad
  \theta=z\frac d{dz}.
$$

Hypergeometric functions arise naturally in many branches of
mathematics. For example, the periods
$$
  \int_1^\infty\frac{dx}{\sqrt{x(x-1)(x-\lambda)}}
$$
of the Legendre family of elliptic curves
$E_\lambda:y^2=x(x-1)(x-\lambda)$ can be expressed as
$$
  {}_2F_1\left(\frac12,\frac12;1;\lambda\right).
$$
Also, it is well-known that
$$
  E_4(\tau)={}_2F_1\left(\frac1{12},\frac5{12};1;\frac{1728}{j(\tau)}\right)^4
$$
where $E_4(\tau)$ is the Eisenstein series of weight $4$ on
$\SL(2,\Z)$ and $j(\tau)$ is the elliptic $j$-function.

In this paper, we are concerned with algebraic transformations of
hypergeometric functions, that is, identities of the form
\begin{equation} \label{equation: general}
  {}_2F_1(a,b;c;z)=R(z){}_2F_1\left(a',b';c';S(z)\right)
\end{equation}
with suitable parameters $a,b,c,a',b',c'$ and algebraic functions
$R(z)$ and $S(z)$. If $w=R(z)$ is of degree $m$ over the field $\C(z)$
or if $z$ is of degree $m$ over the field $\C(w)$, we say the
algebraic transformation has \emph{degree} $m$.
Two examples of algebraic transformations of degree $1$ are given by
$$
  {}_2F_1(a,b;c;z)=(1-z)^{c-a-b}{}_2F_1(c-a,c-b;c;z)
 =(1-z)^{-a}{}_2F_1\left(a,c-b;c;\frac z{z-1}\right).
$$
These identities can be easily proved using the well-known result in
the classical analysis that a Fuchsian differential equation with
precisely $3$ singularities at $0$, $1$, and $\infty$ is completely
determined by the local exponents at these three points.

Beyond transformations of degree $1$, one of the simplest examples is
Kummer's quadratic transformation
$$
  {}_2F_1\left(2a,2b;a+b+\frac12;z\right)
 ={}_2F_1\left(a,b;a+b+\frac12;4z(1-z)\right),
$$
valid for any real numbers $a,b$ with $a+b+1/2\neq 0,-1,-2,\ldots$.
In \cite{Goursat}, Goursat gave more than $100$ algebraic
transformations of degrees $2,3,4,6$. One such example is
$$
  {}_2F_1\left(a,a+\frac13;\frac12;\frac{z(9-8z)^2}{(4z-3)^3}\right)
 =\left(1+\frac z3\right)^{3a}{}_2F_1\left(3a,a+\frac16;\frac12;z\right)
$$
of degree $3$. (See Entry (96) on Page 132 of \cite{Goursat}.) More
recently, Vid\=unas \cite{Vidunas} gave dozens of new algebraic
transformations of degrees $6,8,9,10,12$. For example, he showed that
if we set $\beta=\pm\sqrt{-2}$,
\begin{equation} \label{equation: Vidunas S}
  S(z)=\frac{4z(z-1)(8\beta z+7-4\beta)^8}
  {(2048\beta z^3-3072\beta z^2-3264z^2+912\beta z+3264z+56\beta-17)^3},
\end{equation}
and
$$
  R(z)=\left(1+\frac{16}9(4-17\beta)z-\frac{64}{243}
    (167-136\beta)z^2+\frac{2048}{6561}(112-17\beta)z^3\right)^{-1/16},
$$
then
$$
  {}_2F_1\left(\frac5{24},\frac{13}{24};\frac78;z\right)
 =R(z){}_2F_1\left(\frac1{48},\frac{17}{48};\frac78;S(z)\right),
$$
which is a transformation of degree $10$. (See (32) of
\cite{Vidunas}.) Vid\=unas' examples usually involve Gr\"obner-basis
computation. This is perhaps one of the reasons why Goursat could not
find these transformations.

In this paper, we will present several new algebraic transformations. For
example, one of our favorite identities is
\begin{equation} \label{equation: favorite}
\begin{split}
  &{}_2F_1\left(\frac1{20},\frac14;\frac45;\frac{64z(1-z-z^2)^5}
     {(1-z^2)(1+4z-z^2)^5}\right) \\
  &\qquad\qquad=(1-z^2)^{1/20}(1+4z-z^2)^{1/4}
  {}_2F_1\left(\frac3{10},\frac25;\frac9{10};z^2\right).
\end{split}
\end{equation}
The main novelty in our approach is the interpretation of
${}_2F_1$-hypergeometric functions as automorphic forms on Shimura
curves. Then proving identities such as the one above amounts to
showing two certain automorphic forms on two Shimura curves are
equal. This point of view is especially useful in determining the
function $R(z)$ in \eqref{equation: general}. We will review the basic
definitions regarding Shimura curves and their automorphic forms in
the next section.
\end{section}

\begin{section}{Preliminaries}

The materials in this section are mostly taken from \cite{Vigneras}.
\medskip

\paragraph{\bf Quaternion algebras} A \emph{quaternion algebra} $B$ over a
field $K$ is a central simple algebra of dimension $4$ over $K$. (Here
\emph{central} means that the center of the algebra is $K$ and
\emph{simple} means that $B$ has no proper nontrivial two-sided
ideals.) If the characteristic of $K$ is not $2$, then one can show
that there are elements $i$ and $j$ in $B$ and $a,b\in K^\times$ such
that
$$
  i^2=a, \qquad j^2=b, \qquad ij=-ji,
$$
and $B=K+Ki+Kj+Kij$. In this case, we denote this algebra by
$\JS{a,b}K$. For example, we have $M(2,K)\simeq\JS{1,1}K$ and
$\JS{-1,-1}\R$ is the set of Hamilton's quaternions. Moreover, for
$\alpha=a_0+a_1i+a_2j+a_3ij\in B$, we set
$\overline\alpha=a_0-a_1i-a_2j-a_3ij$. Then the \emph{reduced
  trace} $\tr(\alpha)$ is defined to be
$\alpha+\overline\alpha=2a_0\in K$ and the \emph{reduced norm}
$\nm(\alpha)$ is defined to be
$\alpha\overline\alpha=a_0^2-a_1^2a-a_2^2b+a_3^2ab\in K$.

If $K=\C$, then up to isomorphisms, there is only one quaternion
algebra over $\C$, which is $M(2,\C)$. If $K=\R$ or a non-Archimedean
local field, then up to isomorphism, there are only two quaternion
algebras. One is $M(2,K)$ and the other is a division algebra.

Now assume that $K$ is a number field. Let $v$ be a place of $K$ and
$K_v$ be the completion of $K$ with respect to $v$. If the
localization $B\otimes_KK_v$ is isomorphic to $M(2,K_v)$, we say $B$
\emph{splits} at $v$. If $B\otimes_KK_v$ is isomorphic to a division
algebra, we say $B$ \emph{ramifies} at $v$. It is known that the
number of ramified places is finite and in fact an even integer. The
product of finite ramified places is called the \emph{discriminant} of
the quaternion algebra.

Still assume that $K$ is a number field. Let $R$ be its ring of
integers. An \emph{order} in $B$ is a finite generated $R$-module that
is also a ring with unity containing a basis for $B$. An order is
\emph{maximal} if it is not properly contained in another order. An
\emph{Eichler order} is the intersection of two maximal orders. For
example, if $B=M(2,\Q)$, then $M(2,\Z)$ is a maximal order and
$\SM\Z\Z{N\Z}\Z=M(2,\Z)\cap\SM100N M(2,\Z)\SM100N^{-1}$ is an Eichler
order.
\medskip

\paragraph{\bf Shimura curves} To define a Shimura curve, we assume
that $K$ is a totally real number field and take a quaternion algebra
$B$ over $K$ that splits at exactly one infinite place, that is,
$$
  B\otimes_\Q R\simeq M(2,\R)\times \HH^{[K:\Q]-1},
$$
where $\HH$ is Hamilton's quaternion algebra $\JS{-1,-1}\R$. Then, up
to conjugation, there is a unique embedding $\iota$ of $B$ into
$M(2,\R)$. Note that we have $\nm(\alpha)=\det\iota(\alpha)$ for all
$\alpha\in B$.

Let $\O$ be an order and $\O_1^\ast=\{\alpha\in\O:\nm(\alpha)=1\}$ be
the norm-one group of $\O$. Then the image $\Gamma(\O)$ of $\O_1^\ast$
under the embedding $\iota$ is a discrete subgroup of $\SL(2,\R)$.
Let $\Gamma(\O)$ act on the upper half-plane $\H$ in the usual manner
$$
  \M abcd:\tau\longmapsto\frac{a\tau+b}{c\tau+d}.
$$
Then the quotient space $\Gamma(\O)\backslash\H$ is called the
\emph{Shimura curve} associated to $\O$. For example, if $B=M(2,\Q)$
and $\O=M(2,\Z)$, then $\Gamma(\O)=\SL(2,\Z)$ and
$\Gamma(\O)\backslash\H$ is just the usual modular curve $Y_0(1)$.
Thus, Shimura curves are generalizations of classical modular curves
and they are moduli spaces of certain abelian varieties with
quaternionic multiplication \cite{Shimura-CM1}. In a broader setting,
if $\Gamma$ is any discrete subgroup of $\SL(2,\R)$ commensurable with
$\Gamma(\O)$, then the quotient space $\Gamma\backslash\H$ will also
be called a Shimura curve.

An element $\gamma$ of $\Gamma(\O)$ is \emph{parabolic},
\emph{elliptic}, or \emph{hyperbolic}, according to whether 
$|\tr(\gamma)|=2$, $|\tr(\gamma)|<2$, or $|\tr(\gamma)|>2$. The fixed
point of a parabolic element is called a \emph{cusp}. This can appear
only when $B=M(2,\Q)$. The fixed point $\tau$ of an elliptic element in
$\H$ is called an \emph{elliptic point} of order $n$, where $n$ is the
order of the isotropy subgroup of $\tau$ in $\Gamma(\O)/\pm 1$.

Note that if $B\neq M(2,\Q)$, then the quotient space
$\Gamma(\O)\backslash\H$ is a compact Riemann surface, which we denote
by $X(\O)$. If $B=M(2,\Q)$, we compactify the Riemann surface
$\Gamma(\O)\backslash\H$ by adding cusps and the resulting compact
surface will also be denoted by $X(\O)$.

Now suppose that the compact Riemann surface $X(\O)$ has genus $g$.
Then a classical result says that there exist hyperbolic
elements $A_1,\ldots,A_g$, $B_1,\ldots,B_g$, and elliptic or parabolic
elements $C_1,\ldots,C_r$ that generate $\Gamma(\O)/\pm 1$ with the
single relation
$$
  [A_1,B_1]\ldots[A_g,B_g]C_1\ldots C_r=\mathrm{Id},
$$
where $[A_i,B_i]=A_iB_iA_i^{-1}B_i^{-1}$ is the commutator of $A_i$
and $B_i$. (See \cite[Chapter 4]{Katok}.)
We let $(g;e_1,\ldots,e_r)$ be the \emph{signature} of $X(\O)$.
\medskip

\paragraph{\bf Triangle groups} Suppose that a Shimura curve $X(\O)$
has signature $(0;e_1,e_2,e_3)$. Then we say $\Gamma(\O)$ is an
\emph{arithmetic triangle group}. The complete lists of all arithmetic
triangle groups and their commensurability classes were determined by
Takeuchi \cite{Takeuchi1,Takeuchi2}.

If we cut each fundamental domain of an arithmetic triangle group
$\Gamma(\O)$ into $2$ halves in a suitable way, then the fundamental
half-domains give a tessellation of the upper half-plane $\H$ by
congruent triangles with internal angles $\pi/e_1$, $\pi/e_2$, and
$\pi/e_3$. The following figure shows the tessellation of the unit
disc, which is conformally equivalent to $\H$ through
$\tau\to(\tau-i)/(\tau+i)$, by fundamental half-domains of the
arithmetic triangle group $(0;2,3,7)$.

\centerline{\epsfig{file=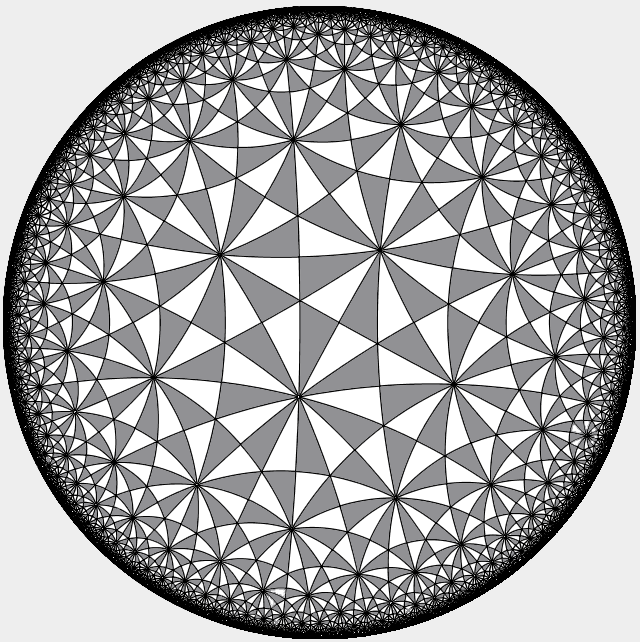,height=2in,width=2in}}

\noindent Here each triangle represents a fundamental half-domain. Any
  combination of a grey triangle with a neighboring white triangle
  will be a fundamental domain for the triangle group $(0;2,3,7)$.
  The triangle group $(0;2,3,7)$ and its associated Shimura curve have
  been studied in details in \cite{Elkies-237}.

In general, for any discrete subgroup $\Gamma$ of $\SL(\R)$ such that
$\Gamma\backslash\H$ has finite volume, we can define its signature in
the same way. If the signature is $(0;e_1,e_2,e_3)$, then we say
$\Gamma$ is a (hyperbolic) triangle group. (There are also notions of
parabolic and elliptic triangle groups, corresponding to tessellation
of $\C$ and $\mathbb P^1(\C)$, respectively.)
\medskip

\paragraph{\bf Automorphic forms on Shimura curves} The definition of
an automorphic form on Shimura curves is the same as that of a modular
form on classical modular curves.

For simplicity, we assume that $B\neq M(2,\Q)$ so that we do not need
to consider cusps. Then an \emph{automorphic form} of weight $k$ on
$\Gamma(\O)$ is a holomorphic function $f:\H\to\C$ such that
\begin{equation} \label{equation: automorphy}
  f\left(\frac{a\tau+b}{c\tau+d}\right)=(c\tau+d)^kf(\tau)
\end{equation}
for all $\SM abcd\in\Gamma(\O)$ and all $\tau\in\H$. The space of
automorphic forms of weight $k$ on $\Gamma(\O)$ will be denoted by
$S_k(\O)$. Also, if a meromorphic function $f:\H\to\C$ satisfies
\eqref{equation: automorphy} with $k=0$, we say $f$ is an
\emph{automorphic function}. If $X(\O)$ has genus $0$, we call an
automorphic function a \emph{Hauptmodul} if it generates the field of
automorphic functions on $\Gamma(\O)$.

Using the Riemann-Roch formula, one can calculate the dimension of
$S_k(\O)$.

\begin{Proposition}[{\cite[Theorem 2.23]{Shimura-book}}]
  \label{proposition: dimension}
  Assume that $B\neq M(2,\Q)$. Suppose that the Shimura curve $X(\O)$
  associated to an order $\O$ in $B$ has signature
  $(g;e_1,\ldots,e_r)$. Then for even integers $k$, we have
  $$
    \dim S_k(\O)=\begin{cases}
    0, &\text{if }k<0, \\
    1, &\text{if }k=0, \\
    g, &\text{if }k=2, \\ \displaystyle
    (k-1)(g-1)+\sum_{j=1}^r\left\lfloor\frac k2\left(1-\frac1{e_j}\right)
    \right\rfloor, &\text{if }k\ge 4. \end{cases}
  $$
\end{Proposition}

The dimension formula for the case $B=M(2,\Q)$ is slightly different.

In the case $B=M(2,\Q)$, there are many methods to construct modular
forms, such as Eisenstein series, theta series, the Dedekind eta
function, and etc. In practice, most explicit methods for modular
curves rely on the Fourier expansions of modular forms and
modular functions, i.e., the expansions with respect to the local
parameter at the cusp $\infty$. However, in the case $B\neq M(2,\Q)$,
because of the lack of cusps on Shimura curves, very few explicit
methods are available for Shimura curves. One of the few methods is
given by the second author of the present paper.

One of the key ideas in \cite{Yang-Schwarzian} is the following
characterization of $S_k(\O)$. Here we assume that the quaternion
algebra is not $M(2,\Q)$.

\begin{Proposition}[{\cite[Theorem 4, Propositions 1 and 6]{Yang-Schwarzian}}]
 \label{theorem: basis}
  Assume that a Shimura curve $X$ has genus zero with elliptic
  points $\tau_1,\ldots,\tau_r$ of order $e_1,\ldots,e_r$,
  respectively. Let $t(\tau)$ be a Hauptmodul of $X$ and set
  $a_i=t(\tau_i)$, $i=1,\ldots,r$. For a positive even integer $k\ge 4$,
  let
  $$
    d_k=\dim S_k(\O)=1-k+\sum_{j=1}^r\left\lfloor
    \frac k2\left(1-\frac1{e_j}\right)\right\rfloor.
  $$
  Then a basis for the space of automorphic forms of weight $k$ on $X$ is
  $$
    t'(\tau)^{k/2}t(\tau)^j\prod_{i=1,a_i\neq\infty}^r
    \left(t(\tau)-a_j\right)^{-\lfloor k(1-1/e_1)/2\rfloor}, \quad
    j=0,\ldots,d_k-1.
  $$
  Moreover, the functions $t'(\tau)^{1/2}$ and $\tau t'(\tau)^{1/2}$,
  as functions of $t$, satisfy the differential equation
  $$
    f''+Q(t)f=0,
  $$
  where
  $$
   Q(t)=\frac14\sum_{j=1,a_j\neq\infty}^r\frac{1-1/e_j^2}{(t-a_j)^2}
        +\sum_{j=1,a_j\neq\infty}^r\frac{B_j}{t-a_j}
  $$
  for some constants $B_j$. Moreover, if $a_j\neq\infty$ for all $j$,
  then the constants $B_j$ satisfy
  $$
  \sum_{j=1}^r B_j=
  \sum_{j=1}^r\left(a_jB_j+\frac14(1-1/e_j^2)\right)=
  \sum_{j=1}^r\left(a_j^2B_j+\frac12a_j(1-1/e_j^2)\right)=0.
  $$
  Also, if $a_r=\infty$, then $B_j$ satisfy
  $$
  \sum_{j=1}^{r-1}B_j=0, \qquad
  \sum_{j=1}^{r-1}\left(a_jB_j+\frac14(1-1/e_j^2)\right)=\frac14(1-1/e_r^2).
  $$
\end{Proposition}

\begin{Remark} In \cite{Yang-Schwarzian}, the differential equation
  $f''+Q(t)f=0$ is called the \emph{Schwarzian differential equation}
  associated to $t$ because $Q(t)$ is related to the Schwarzian
  derivative by the relation
  $$
    2Q(t)t'(\tau)^2+\{t,\tau\}=0,
  $$
  where
  $$
    \{t,\tau\}=\frac{t'''(\tau)}{t'(\tau)}-\frac32\left(
    \frac{t''(\tau)}{t'(\tau)}\right)^2
  $$
  is the \emph{Schwarzian derivative}. In general, in literature
  \cite{Bayer-Travesa}, if $f$ is a thrice-differentiable function of
  $z$, then
  $$
    D(f,z):=-\frac{\{f,z\}}{2f'(z)^2}
  $$
  is called the \emph{automorphic derivative} associated to $f$ and
  $z$. In the case $f$ is an automorphic function on a Shimura curve,
  then $D(f,\tau)$ is also an automorphic function. In particular, if
  $t$ is a Hauptmodul on a Shimura curve of genus $0$, then
  $Q(t)=D(t,\tau)$ is a rational function of $t$.
\end{Remark}

The upshot of this result is that it is often possible to determine
the differential equation without explicitly constructing a
Hauptmodul. For example, if $\Gamma(\O)$ is a triangle group with
signature $(0;e_1,e_2,e_3)$, then there always exists a (unique)
Hauptmodul $t$ with $a_1=0$, $a_2=1$, and $a_3=\infty$. Then the
relations between $B_j$ uniquely determine the differential equation.
In general, one can usually use coverings between Shimura curves of
genus $0$ to determine the differential equation. This is done by the
first author \cite{Tu} of the present paper for many Shimura curves of
genus $0$ associated to orders in quaternion algebras over $\Q$. Once
the differential equation is determined, one can express automorphic
forms in terms of $t$-series and then study properties of automorphic
forms using these $t$-series. For example, in \cite{Yang-Schwarzian} the
second author devised a method to compute Hecke operators on these
$t$-series.

In the case of triangle groups, since the number of singularities of
the differential equation is $3$, the differential equation is
essentially a hypergeometric differential equation.

\begin{Proposition}[{\cite[Theorem 9]{Yang-Schwarzian}}]
  \label{theorem: triangle}
  Assume that a Shimura curve $X$ has signature
  $(0;e_1,e_2,e_3)$. Let $t(\tau)$ be the Hauptmodul of $X$ with
  values $0$, $1$, and $\infty$ at the elliptic points of order $e_1$,
  $e_2$, and $e_3$, respectively. Let $k\ge 4$ be an even integer.
  Then a basis for the space of automorphic forms of weight $k$ on $X$
  is given by
  $$
    t^{\{k(1-1/e_1)/2\}}(1-t)^{\{k(1-1/e_2)/2\}}t^j\left(
    {}_2F_1(a,b;c;t)+Ct^{1/e_1}{}_2F_1(a',b',c';t)\right)^k,
  $$
  $j=0,\ldots,\lfloor k(1-1/e_1)/2\rfloor+\lfloor k(1-1/e_2)/2\rfloor
  +\lfloor k(1-1/e_3)/2\rfloor-k$, for some constant $C$, where
  for a rational number $x$, we let $\{x\}$ denote the fractional part
  of $x$,
  $$
    a=\frac12\left(1-\frac1{e_1}-\frac1{e_2}-\frac1{e_3}\right), \qquad
    b=a+\frac1{e_3}, \qquad c=1-\frac1{e_1}
  $$
  and
  $$
    a'=a+\frac1{e_1}, \qquad b'=b+\frac1{e_1}, \qquad
    c'=c+\frac2{e_1}.
  $$
\end{Proposition}

For general Shimura curves of genus $0$, the following properties of
automorphic derivatives are very useful in determining the Schwarzian
differential equation associated to a Hauptmodul.

\begin{Proposition} \label{proposition: D}
  Automorphic derivatives have the following properties.
  \begin{enumerate}
  \item $D((az+b)/(cz+d),z)=0$ for all $\SM abcd\in\mathrm{GL}(2,\C)$.
  \item $D(g\circ f,z)=D(g,f(z))+D(f,z)/(dg/df)^2$.
 \end{enumerate}
\end{Proposition}

\begin{Proposition} \label{proposition: D 1}
  Let $z(\tau)$ be a Hauptmodul for a Shimura curve $X(\O)$ of genus
  $0$. Let $R(x)\in\C(x)$ be the rational function such that the
  automorphic derivative $Q(z)=D(z,\tau)$ is equal to $R(z)$.
  Assume that $\gamma$ is an element of $\SL(2,\R)$ normalizing the
  norm-one group of $\O$ and let $\sigma$ be the automorphism of
  $X(\O)$ induced by $\gamma$. If $\sigma:z\mapsto(az+b)/(cz+d)$, then
  $R(x)$ satisfies
  $$
    \frac{(ad-bc)^2}{(cx+d)^4}R\left(\frac{ax+b}{cx+d}\right)=R(x).
  $$ 
\end{Proposition}

\begin{proof} We shall compute $D(z(\gamma\tau),\tau)$ in two ways.
  By Proposition \ref{proposition: D}, we have
  \begin{equation*}
  \begin{split}
    D(z(\gamma\tau),\tau)=D\left(\frac{az(\tau)+b}{cz(\tau)+d},z(\tau)\right)
      +\frac{D(z(\tau),\tau)}{(dz(\gamma\tau)/dz(\tau))^2}
    =0+\frac{(cz+d)^4R(z)}{(ad-bc)^2}.
  \end{split}
  \end{equation*}
  On the other hand, by the same proposition, we also have
  \begin{equation*}
  \begin{split}
    D(z(\gamma\tau),\tau)=D(z(\gamma\tau),\gamma\tau)
     +\frac{D(\gamma\tau,\tau)}{(dz(\gamma\tau)/d\gamma\tau)^2}
    =R(z(\gamma\tau))=R\left(\frac{az+b}{cz+d}\right).
  \end{split}
  \end{equation*}
  Comparing the two expressions, we get the formula.
\end{proof}
\medskip

\paragraph{\bf Algebraic transformations of hypergeometric functions}
Consider the following situation. Suppose that $\Gamma_1<\Gamma_2$ are
two arithmetic triangle groups with Hauptmoduls $z_1$ and $z_2$,
respectively. Since any automorphic function on $\Gamma_2$ is also an
automorphic function on $\Gamma_1$, we have $z_2=S(z_1)$ for some
$S(x)\in\C(x)$. Likewise, if $f_1$ and $f_2$ are two automorphic forms
of the same weight $k$ on $\Gamma_1$ and $\Gamma_2$, respectively, then
the ratio $f_1/f_2$ is an automorphic function on $\Gamma_1$ and hence
is equal to $R(z_1)$ for some $R(x)\in\C(x)$. In view of Proposition
\ref{theorem: triangle}, after taking the $k$th roots of the two sides
of $f_1/f_2=R(z_1)$, we obtain an algebraic transformation of
hypergeometric function. This explains the existence of Kummer's,
Goursat's and Vid\=unas' transformations. (Of course, the triangle
groups appearing in their transformations may not be arithmetic, but
the argument above is still valid.)

More generally, if $\Gamma_1$ and $\Gamma_2$ are two commensurable
arithmetic triangle groups such that the Shimura curve associated to
$\Gamma=\Gamma_1\cap\Gamma_2$ has genus $0$. Let $z$ be a Hauptmodul
on $\Gamma$. Then each of $z_1$ and $z_2$ is a rational function of
$z$. Similarly, the ratio $f_1/f_2$ is also a rational function of
$z$. Again, Proposition \ref{theorem: triangle} yields an algebraic
transformation of the form
$$
  _2F_1\left(a_1,b_1;c_1;S_1(z)\right)=R(z)
  {}_2F_1\left(a_2,b_2;c_2;S_2(z)\right)
$$
for some rational functions $S_1(z)$ and $S_2(z)$ and some algebraic
function $R(z)$. This is the theory behind \eqref{equation: favorite}
and other algebraic transformations given in the paper.

\begin{Definition} Let $S(z)\in\C(z)$ be a rational function. If the
  finite covering $\P^1(\C)\to\P^1(\C)$ defined by $S:z\to S(z)$ is
  ramified at most at three points $0$, $1$, and $\infty$, then $S$ is
  called a \emph{Belyi function}.
\end{Definition}

In practice, the Belyi functions $S_1(z)$ and $S_2(z)$ can be
determined by the ramification data of the coverings of Shimura
curves. The function $R(z)$ can be determined by Propositions
\ref{theorem: basis} and \ref{theorem: triangle}.

We now obtain algebraic transformations of hypergeometric functions
using the above idea. Note that according to \cite{Takeuchi2},
arithmetic triangle groups fall in $19$ commensurability classes.
The first class in his list corresponds to classical modular curves.
In this case, it is easier to use classical modular forms to derive
identities. We will not discuss this case here. Identities arising
from Classes II, V, and XII are special cases of a family of
identities, and so are identities from Classes IV, VIII, XI, XIII, XV,
and XVII. These cases will be treated in a later section. Here we
first consider Class III in Section \ref{section: class III} and Class
VI in Section \ref{section: class VI}. (There are no identities from
Classes IX and XIX since these classes consist of a single group.
Also, identities from Classes VII, XIV, XVI, and XVIII are just
Kummer's quadratic transformations.)
\end{section}

\begin{section}{Algebraic transformations associated to Class III}
\label{section: class III}

According to \cite{Takeuchi2}, Takeuchi's Class III of commensurable
arithmetic triangle groups has the following subgroup diagram. Here
because all groups involved have genus zero, we omit the genus
information in the signatures of the groups.
$$
  \begin{diagram}
  \node{(2,6,8)} \arrow{s,l,-}{2} \arrow{ese,r,-}{2}
  \node[4]{(2,3,8)} \arrow{wsw,l,-}{10} \arrow{s,l,-}{2}
    \arrow{ese,r,-}{3} \\
  \node{(4,6,6)}
  \node[2]{(3,8,8)}
  \node[2]{(3,3,4)} \arrow{s,r,-}{3}
  \node[2]{(2,4,8)} \arrow{wsw,r,-}{2} \arrow{s,r,-}{2} \\
  \node[5]{(4,4,4)}
  \node[2]{(2,8,8)} \arrow{s,r,-}{2} \\
  \node[7]{(4,8,8)}
 \end{diagram}
$$
The main goal in this section is to prove an algebraic transformation
associated to the pair of triangle groups $(4,6,6)$ and $(4,4,4)$.

\begin{Theorem} \label{theorem: Class III} Let $\alpha$ be a root of
  $x^2+3=0$ and $\beta$ a root of $x^2+2=0$. We have
\begin{equation}
\begin{split}
 &\frac{(1+z)^{1/8}(1-3z)^{1/8}}{(1+\alpha z)^{5/4}}
  {}_2F_1\left(\frac5{24},\frac38;\frac34;
  \frac{12\alpha z(1-z^2)(1-9z^2)}{(1+\alpha z)^6}\right) \\
 &\qquad=\frac{1}{(1+(4+2\beta)z-(1+2\beta)z^2)^{1/2}}
  {}_2F_1\left(\frac18,\frac38;\frac34;R(z)\right), 
\end{split}
\end{equation}
and
\begin{equation}
\begin{split}
  &\frac{(1-z)^{1/4}(1+z)^{5/8}(1-3z)^{1/4}(1+3z)^{5/8}}{(1+\alpha z)^{11/4}}
   {}_2F_1\left(\frac{11}{24},\frac58;\frac54;
   \frac{12\alpha z(1-z^2)(1-9z^2)}{(1+\alpha z)^6}\right) \\
  &\qquad=\frac{(1+(-7+4\beta)z^2/3)}
   {(1+(4+2\beta)z-(1+2\beta)z^2)^{3/2}}
   {}_2F_1\left(\frac38,\frac58;\frac54;R(z)\right)
\end{split}
\end{equation}
where
$$
  R(z)=-\frac{4(1+\beta)^4z(1+(-7+4\beta)z^2/3)^4}
  {(1+z)(1-3z)(1+(4+2\beta)z-(1+2\beta)z^2)^4}.
$$
\end{Theorem}

We first determine the signatures of the intersections.

\begin{Lemma} \label{lemma: Class III subgroups} We have
$$
  \begin{diagram}
  \node{(2,6,8)} \arrow{s,l,-}{2} \arrow{ese,r,-}{2}
  \node[4]{(2,3,8)} \arrow{wsw,l,-}{10} \arrow{s,l,-}{2}
    \arrow{ese,r,-}{3^\ast} \\
  \node{\Gamma_1=(4,6,6)} \arrow{ese,l,-}{2}
  \node[2]{\Gamma_2=(3,8,8)} \arrow{s,l,-}{2}
  \node[2]{\Gamma_3=(3,3,4)} \arrow{wsw,l,-}{10} \arrow{s,r,-}{3}
  \node[2]{(2,4,8)} \arrow{wsw,r,-}{2} \arrow{s,r,-}{2} \\
  \node[3]{\Gamma_5=(3,4,3,4)} \arrow{s,l,-}{3}
  \node[2]{\Gamma_4=(4,4,4)} \arrow{wsw,l,-}{10}
  \node[2]{(2,8,8)} \arrow{s,r,-}{2} \\
  \node[3]{\Gamma_6=(4^6)}
  \node[4]{(4,8,8)}
  \end{diagram}
$$
Moreover, the group of signature $(4^6)$ is a normal subgroup of the
group of signature $(3,4,3,4)$. (Here $(4^6)$ is a shorthand for
$(4,4,4,4,4,4)$.)
\end{Lemma}

\begin{proof} Let $\Gamma_2=(3,8,8)$ and $\Gamma_2'$ be its commutator
  subgroup. From the group presentation
  $$
    \Gamma_2\simeq\gen{\gamma_1,\gamma_2:\gamma_1^3,\gamma_2^8,(\gamma_1\gamma_2)^8}
  $$
  for $\Gamma_2$, we know that $\Gamma_2/\Gamma_2'$ is cyclic of order
  $8$. Thus, $\Gamma_2$ has exactly one subgroup of index $2$, which
  must be the common intersection of the groups $(4,6,6)$, $(3,8,8)$
  and $(3,3,4)$. The signature of this subgroup can be easily
  determined by observing that a covering of degree $2$ from a Shimura
  curve to the Shimura curve associated to $(3,8,8)$ can only ramify
  at the two elliptic points of order $8$. We find that the signature must be
  $(3,4,3,4)$.

  We next observe that the commutator subgroup $\Gamma_3'$ of the
  group $\Gamma_3=(3,3,4)$ is cyclic of order $3$. Thus, $\Gamma_3'$
  is a normal subgroup of index $3$ in $\Gamma_3$. This $\Gamma_3'$
  must be the same as the group of signature $(4,4,4)$. If
  $\Gamma_3'\neq(4,4,4)$, then $\Gamma_3'\cap(4,4,4)$ is a normal
  subgroup of $(4,4,4)$ of index $3$, but the group $(4,4,4)$ cannot
  have a normal subgroup of index $3$. We next determine the signature
  of the intersection of $\Gamma_4=(4,4,4)$ and $\Gamma_5=(3,4,3,4)$.

  Let $X_j$ denote the Shimura curve associated to the group $\Gamma_j$.
  Since $\Gamma_4$ is a normal subgroup of $\Gamma_3$ of index $3$,
  the intersection $\Gamma_6$ of $\Gamma_4$ and $\Gamma_5$ is a normal
  subgroup of index $3$ in $\Gamma_5$, which implies that the two
  elliptic points of order $4$ of $X_5$ must split completely on
  $X_6$. In view of the Riemann-Hurwitz formula, the two elliptic
  points of order $3$ of $X_5$ must be totally ramified. We conclude
  that $\Gamma_6$ has signature $(4^6)$.

  In fact, the subgroup relations mentioned above can be visualized by
  the following figures.

\centerline{
  \epsfig{file=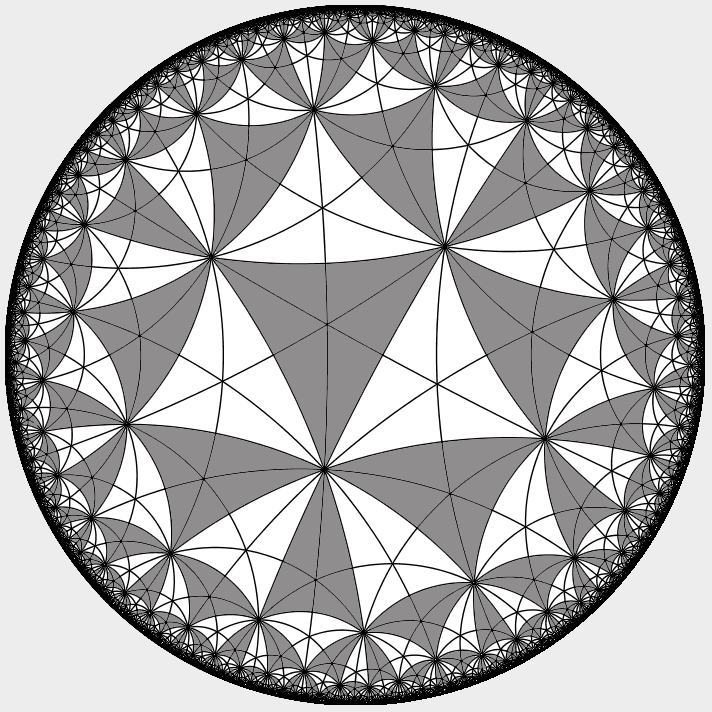,width=1.6in,height=1.6in} \qquad
  \epsfig{file=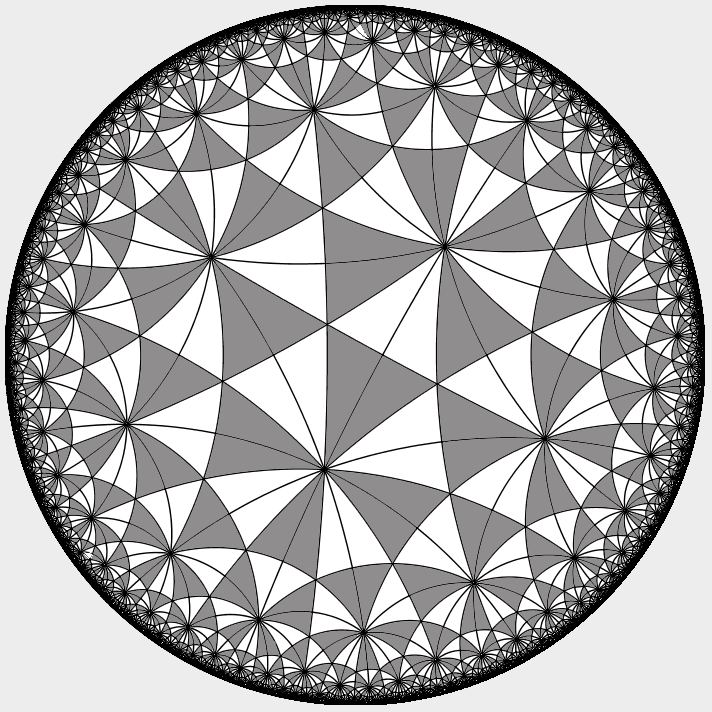,width=1.6in,height=1.6in}
}

\noindent Here the small triangles are $(2,3,8)$-triangles. Let $G$ be
the group of all symmetries of the tessellation of the hyperbolic
plane by the $(4,4,4)$-triangles and $G_0$ be the subgroup generated
by the reflections across the edges of $(4,4,4)$-triangles. Then
$G/G_0$ is isomorphic to $D_3$. The $(3,3,4)$-triangle group
corresponds to the cyclic subgroup of order $3$ in $G/G_0$, while the
group $(2,3,8)$ corresponds the whole group $G/G_0$. Similarly, if we piece $12$
copies of $(2,6,8)$-triangles around the vertex of inner angle
$\pi/4$, we get a regular hexagon with inner angles $\pi/4$. Let $H$
be the group of all symmetries of the tessellation by this regular
hexagon and $H_0$ be the subgroup generated by the reflections across
the edges of hexagons. Then $H/H_0$ is isomorphic to $D_6$. The unique
cyclic subgroup of order $3$ in $H/H_0$ corresponds to the group
$(3,4,3,4)$. See the figures below.

\centerline{
  \epsfig{file=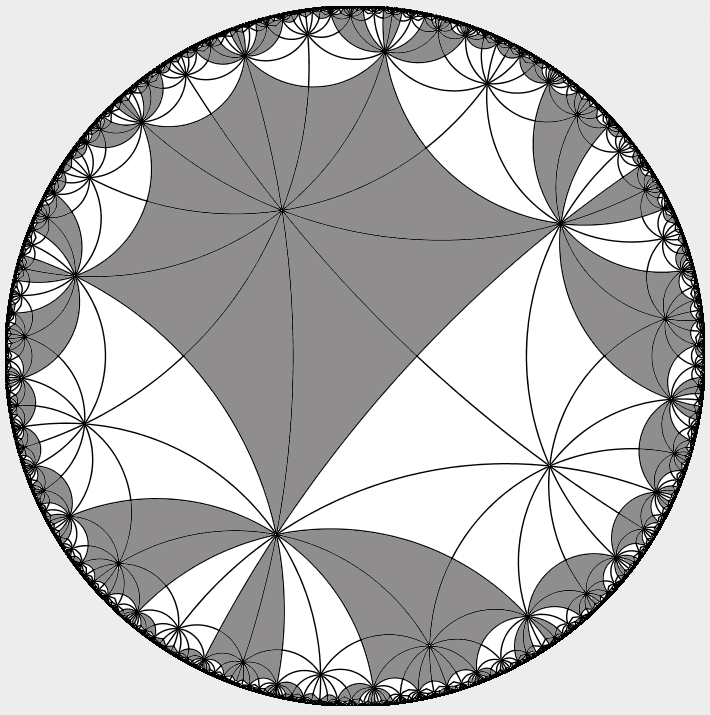,width=1.6in,height=1.6in} \qquad
  \epsfig{file=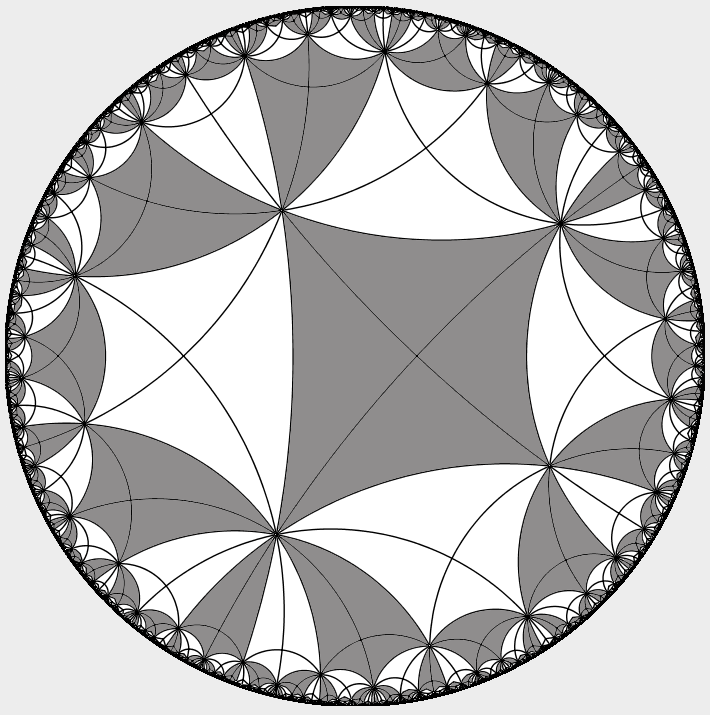,width=1.6in,height=1.6in}
}

\noindent (The groups $(2,6,8)$, $(4,6,6)$, and $(3,8,8)$ correspond
to the whole $H/H_0$, the cyclic subgroup of order $6$ of $H/H_0$, and
one of the $D_3$-subgroups, respectively.)
\end{proof}

Now let $\Gamma_1=(4,6,6)$, $\Gamma_2=(3,8,8)$, $\Gamma_3=(3,3,4)$,
$\Gamma_4=(4,4,4)$, $\Gamma_5=(3,4,3,4)$, and $\Gamma_6=(4^6)$. Let
$X_j=X(\Gamma_j)$, $j=1,\ldots,6$, be the corresponding Shimura
curves. Label the elliptic points on $X_1$ by $P_4$, $P_6$, and
$P_6'$, those on $X_2$ by $Q_3$, $Q_8$, and $Q_8'$, those on $X_3$ by
$R_3$, $R_3'$, and $R_4$, those on $X_4$ by $S_4$, $S_4'$, $S_4''$,
and those on $X_5$ by $T_3$, $T_3'$, $T_4$, and $T_4'$ (with the
subscripts carrying the obvious meaning) such that the ramification
data are given by

\centerline{
\begin{tikzpicture}[scale=0.06]
  \node[place] (T31) [label=above:$T_3$, token=1] at (-75,20) {};
  \node[place] (T32) [label=above:$T_3'$, token=1] at (-55,20) {};
  \node[place] (T41) [label=above:$T_4$, token=1] at (-35,20) {};
  \node[place] (T42) [label=above:$T_4'$, token=1] at (-15,20) {};
  \node[place] (T311) [label=above:$T_3$, token=1] at (15,20) {};
  \node[place] (T321) [label=above:$T_3'$, token=1] at (35,20) {};
  \node[place] (T411) [label=above:$T_4$, token=1] at (55,20) {};
  \node[place] (T421) [label=above:$T_4'$, token=1] at (75,20) {};
  \node[place] (P61) [label=below:$P_6$, token=1] at (-75,0) {};
  \node[place] (P62) [label=below:$P_6'$, token=1] at (-55,0) {};
  \node[place] (P4) [label=below:$P_4$, token=1] at (-25,0) {};
  \node[place] (Q3) [label=below:$Q_3$, token=1] at (25,0) {};
  \node[place] (Q81) [label=below:$Q_8$, token=1] at (55,0) {};
  \node[place] (Q82) [label=below:$Q_8'$, token=1] at (75,0) {};
  \draw[line width=.6pt] (T31) --  node[left=2pt] {$2$} (P61);
  \draw[line width=.6pt] (T32) --  node[left=2pt] {$2$}(P62);
  \draw[line width=.6pt] (T41) --  node[left=2pt] {$1$}(P4);
  \draw[line width=.6pt] (T42) --  node[right=2pt] {$1$}(P4);
  \draw[line width=.6pt] (T311) -- node[left=2pt] {$1$}(Q3);
  \draw[line width=.6pt] (T321) -- node[right=2pt] {$1$}(Q3);
  \draw[line width=.6pt] (T411) -- node[right=2pt] {$2$}(Q81);
  \draw[line width=.6pt] (T421) -- node[right=2pt] {$2$}(Q82);
\end{tikzpicture}
}
\centerline{
\begin{tikzpicture}[scale=0.05]
  \node[place] (T31) [label=above:$T_3$, token=1] at (-110,20) {};
  \node[place] (A1) [token=1] at (-90,20) {};
  \node[place] (A2) [token=1] at (-70,20) {};
  \node[place] (A3) [token=1] at (-50,20) {};
  \node[place] (T32) [label=above:$T_3'$, token=1] at (-30,20) {};
  \node[place] (B1) [token=1] at (-10,20) {};
  \node[place] (B2) [token=1] at (10,20) {};
  \node[place] (B3) [token=1] at (30,20) {};
  \node[place] (T41) [label=above:$T_4$, token=1] at (50,20) {};
  \node[place] (T42) [label=above:$T_4'$, token=1] at (70,20) {};
  \node[place] (C1) [token=1] at (90,20) {};
  \node[place] (C2) [token=1] at (110,20) {};
  \node[place] (R31) [label=below:$R_3$, token=1] at (-80,0) {};
  \node[place] (R32) [label=below:$R_3'$, token=1] at (0,0) {};
  \node[place] (R4) [label=below:$R_4$, token=1] at (80,0) {};
  \draw[line width=.6pt] (T31) --  node[left=2pt] {$1$} (R31);
  \draw[line width=.6pt] (A1) --  node[left=2pt] {$3$} (R31);
  \draw[line width=.6pt] (A2) --  node[right=2pt] {$3$}(R31);
  \draw[line width=.6pt] (A3) --  node[right=2pt] {$3$}(R31);
  \draw[line width=.6pt] (T32) -- node[left=2pt] {$1$}(R32);
  \draw[line width=.6pt] (B1) -- node[left=2pt] {$3$}(R32);
  \draw[line width=.6pt] (B2) -- node[right=2pt] {$3$}(R32);
  \draw[line width=.6pt] (B3) -- node[right=2pt] {$3$}(R32);
  \draw[line width=.6pt] (T41) -- node[left=2pt] {$1$} (R4);
  \draw[line width=.6pt] (T42) -- node[left=2pt] {$1$} (R4);
  \draw[line width=.6pt] (C1) -- node[right=2pt] {$4$} (R4);
  \draw[line width=.6pt] (C2) -- node[right=2pt] {$4$} (R4);
\end{tikzpicture}
}
Label the elliptic points of $X_6$ by $U_1,\ldots,U_6$ such that the
rotation around the center of the $(4^6)$-polygon by the angle $\pi/3$
permutes the six points cyclically. From the figures above, we know
that if we label the points such that $U_1$ lies above $T_4$, then
the ramification data for $X_6\to X_5$ are

\centerline{
\begin{tikzpicture}[scale=0.06]
  \node[place] (U1) [label=above:$U_1$, token=1] at (-70,20) {};
  \node[place] (U3) [label=above:$U_3$, token=1] at (-50,20) {};
  \node[place] (U5) [label=above:$U_5$, token=1] at (-30,20) {};
  \node[place] (U2) [label=above:$U_2$, token=1] at (-10,20) {};
  \node[place] (U4) [label=above:$U_4$, token=1] at (10,20) {};
  \node[place] (U6) [label=above:$U_6$, token=1] at (30,20) {};
  \node[place] (U01) [label=above:$U_0$, token=1] at (50,20) {};
  \node[place] (U02) [label=above:$U_0'$, token=1] at (70,20) {};
  \node[place] (T31) [label=below:$T_4$, token=1] at (-50,0) {};
  \node[place] (T32) [label=below:$T_4'$, token=1] at (10,0) {};
  \node[place] (T41) [label=below:$T_3$, token=1] at (50,0) {};
  \node[place] (T42) [label=below:$T_3'$, token=1] at (70,0) {};
  \draw[line width=.6pt] (U1) --  node[left=2pt] {$1$} (T31);
  \draw[line width=.6pt] (U3) --  node[left=1pt] {$1$} (T31);
  \draw[line width=.6pt] (U5) --  node[right=2pt] {$1$} (T31);
  \draw[line width=.6pt] (U2) --  node[left=2pt] {$1$} (T32);
  \draw[line width=.6pt] (U4) --  node[left=1pt] {$1$} (T32);
  \draw[line width=.6pt] (U6) --  node[right=2pt] {$1$} (T32);
  \draw[line width=.6pt] (U01) -- node[right=2pt] {$3$} (T41);
  \draw[line width=.6pt] (U02) -- node[right=2pt] {$3$} (T42);
\end{tikzpicture}
}
\noindent where $U_0$ and $U_0'$ are the centers of the
$(4^6)$-polygons. (The reader is reminded that each $(4^6)$-polygon
represents only half of the fundamental domain for the Shimura curve
$X_6$. Referring to the figure in the proof of the lemma above, a
fundamental domain consists of a grey $(4^6)$-polygon and a
neighboring white $(4^6)$-polygon.)

\begin{Lemma} \label{lemma: Class III Belyi}
  The two elliptic points of $X_6$ at the two ends of a
  diagonal of a $(4^6)$-polygon lie above the same elliptic point of
  $X_4$. That is, labeling the elliptic points of $X_4$ suitably,
  we have

\centerline{
\begin{tikzpicture}[scale=0.05]
  \node[place] (U1) [label=above:$U_1$, token=1] at (-110,20) {};
  \node[place] (U4) [label=above:$U_4$, token=1] at (-90,20) {};
  \node[place] (A1) [token=1] at (-70,20) {};
  \node[place] (A2) [token=1] at (-50,20) {};
  \node[place] (U2) [label=above:$U_2$, token=1] at (-30,20) {};
  \node[place] (U5) [label=above:$U_5$, token=1] at (-10,20) {};
  \node[place] (B1) [token=1] at (10,20) {};
  \node[place] (B2) [token=1] at (30,20) {};
  \node[place] (U3) [label=above:$U_3$, token=1] at (50,20) {};
  \node[place] (U6) [label=above:$U_6$, token=1] at (70,20) {};
  \node[place] (C1) [token=1] at (90,20) {};
  \node[place] (C2) [token=1] at (110,20) {};
  \node[place] (S1) [label=below:$S_4$, token=1] at (-80,0) {};
  \node[place] (S2) [label=below:$S_4'$, token=1] at (0,0) {};
  \node[place] (S3) [label=below:$S_4''$, token=1] at (80,0) {};
  \draw[line width=.6pt] (U1) --  node[left=2pt] {$1$} (S1);
  \draw[line width=.6pt] (U4) --  node[left=1pt] {$1$} (S1);
  \draw[line width=.6pt] (A1) --  node[right=1pt] {$4$}(S1);
  \draw[line width=.6pt] (A2) --  node[right=2pt] {$4$}(S1);
  \draw[line width=.6pt] (U2) -- node[left=2pt] {$1$}(S2);
  \draw[line width=.6pt] (U5) -- node[left=1pt] {$1$}(S2);
  \draw[line width=.6pt] (B1) -- node[right=1pt] {$4$}(S2);
  \draw[line width=.6pt] (B2) -- node[right=2pt] {$4$}(S2);
  \draw[line width=.6pt] (U3) -- node[left=2pt] {$1$} (S3);
  \draw[line width=.6pt] (U6) -- node[left=1pt] {$1$} (S3);
  \draw[line width=.6pt] (C1) -- node[right=1pt] {$4$} (S3);
  \draw[line width=.6pt] (C2) -- node[right=2pt] {$4$} (S3);
\end{tikzpicture}
}
Moreover, if we choose Hauptmoduls $z_j(\tau)$ for $X_j$,
$j=1,\ldots,6$, by requiring
\begin{equation*}
\begin{split}
  z_1(P_4)=0, \quad &z_1(P_6)=1, \quad z_1(P_6')=\infty, \\
  z_2(Q_8)=0, \quad &z_2(Q_3)=1, \quad z_2(Q_8')=\infty, \\
  z_3(R_4)=0, \quad &z_3(R_3)=1, \quad z_3(R_3')=\infty, \\
  z_4(S_4)=0, \quad &z_4(S_4')=1, \quad z_4(S_4'')=\infty, \\
  z_5(T_4)=0, \quad &z_5(T_3)=1, \quad z_5(T_4')=\infty, \\
  z_6(U_1)=0, \quad &z_6(U_3)=1, \quad z_6(U_4)=\infty,
\end{split}
\end{equation*}
then we have
\begin{equation*}
\begin{split}
  z_1=\frac{4z_5}{(1+z_5)^2},& \qquad z_2=z_5^2, \\
  z_3=\frac{3(\zeta-\zeta^2)z_4(1-z_4)}{(1+\zeta z_4)^3}, &\qquad
  z_5=\frac{3(\zeta-\zeta^2)z_6(1-z_6^2)}{1-9z_6^2},
\end{split}
\end{equation*}
$$
  z_3=\frac{(28+16\beta)z_5(1+(-17+56\beta)z_5^2/81)^4}
  {(1+z_5)(1+(13+8\beta)z_5/3-(25+32\beta)z_5^2/9
  +(17-56\beta)z_5^3/81)^3},
$$
and
$$
  z_4=-\frac{4(1+\beta)^4z_6(1+(-7+4\beta)z_6^2/3)^4}
  {(1+z_6)(1-3z_6)(1+(4+2\beta)z_6-(1+2\beta)z_6^2)^4},
$$
where $\zeta$ is a $3$rd root of unity and $\beta$ is a root of
$x^2+2=0$.
\end{Lemma}

\begin{proof} The ramification data for the covering $X_5\to X_2$ and
  the assumption $z_2(Q_3)=z_5(T_3)=1$ imply that $z_2=z_5^2$ and
  $$
    z(T_3')=-1.
  $$
  The relation between $z_1$ and $z_5$ is easy to determine. We find
  $z_1=4z_5/(1+z_5)^2$.

  To determine the relation between $z_3$ and $z_4$, we recall from
  Lemma \ref{lemma: Class III subgroups} that $\Gamma_4$ is a normal
  subgroup $\Gamma_3$. Any element of $\Gamma_3$ not in $\Gamma_4$
  induces an automorphism of order $3$ on $X_4$. Such an automorphism
  must permute the three elliptic points $S_4$, $S_4'$, and $S_4''$
  cyclically. In term of the Hauptmodul $z_4$, such an automorphism is
  either
  $$
    \sigma: z_4\longmapsto \frac{-1}{z_4-1}
  $$
  or its square. Moreover, the fixed points of such an automorphism
  are the ramified points in the covering $X_4\to X_3$. That is, if we
  let $S_0$ and $S_0'$ be the points lying above $R_3$ and $R_3'$
  respectively, then $z_4(S_0),z_4(S_0')\in\{-\zeta,-\zeta^2\}$, where
  $\zeta$ is a primitive $3$rd root of unity. Then from the
  ramification data, we easily deduce that
  $z_3=(\zeta-\zeta^2)z_4(1-z_4^2)/(1+\zeta z_4)^3$.

  To determine the relation between $z_5$ and $z_6$, we argue
  similarly as above. The tessellation of the hyperbolic plane by
  $\Gamma_6$ has a $D_6$-symmetry, in addition to the symmetries
  arising from the reflections across the edges of the
  $(4^6)$-polygons. This provides many useful informations. For
  example, if we let $\tau$ be the reflection across the diagonal
  joining $U_1$ and $U_4$, then $\tau$ induces an involution on $X_6$,
  which, in terms of $z_6$, is given by
  $$
    \tau:z_6\longmapsto-z_6,
  $$
  which implies that
  $$
    z_6(P_5)=-1.
  $$
  Furthermore, let $\rho$ denote the rotation by angle $\pi/3$ around
  the center of the hexagon. Then
  $$
    \rho:z_6\longmapsto\frac{cz_6+1}{-cz_6+c}
  $$
  for some zero constant $c$ since $\rho$ maps $1$ to $\infty$ and
  $\infty$ to $-1$. In light of $\rho^2:0\to 1$, we conclude that
  $c=3$ and
  $$
    z_6(U_2)=1/3, \qquad z_6(U_6)=-1/3.
  $$
  It follows that $z_5=Az_6(1-z_6^2)/(1-9z_6^2)$ for some $A$. This
  constant $A$ has the property that $Ax(1-x^2)-(1-9x^2)$ has repeated
  roots. We find $A=\pm 3\sqrt{-3}$. The choice of the sign must be
  synchronized with the choice of the third root of unity in the
  relation between $z_4$ and $z_5$. This will be done later.

  We now come to the more complicated part of the lemma. Let
  $\pi:X_6\to X_4$ be the covering of the Shimura curves. Let $\gamma$
  be an element of $\Gamma_5$ not in $\Gamma_6$. Then $\gamma$
  normalizes both $\Gamma_4$ and $\Gamma_6$ and induces automorphisms
  $\rho_1$ and $\rho_2$ on $X_4$ and $X_6$, respectively. We may
  assume that $\rho_2=\rho^2$, where $\rho$ permutes $U_1,\ldots,U_6$
  cyclically, as defined in the previous paragraph. It is easy
  to check that $\pi\circ\rho_1=\rho_2\circ\pi$. Thus, $\pi(U_1)$,
  $\pi(U_3)$, and $\pi(U_5)$ are three different elliptic points on
  $X_4$. We label them by $S_4$, $S_4'$, and $S_4''$, respectively.
  Let $V_1,V_2$ be the two ramified points lying above $S_4$.
  Now there are three possibilities
  $$
    \pi^{-1}(S_4)=\{U_1,U_2,V_1,V_2\}, \,
    \pi^{-1}(S_4)=\{U_1,U_4,V_1,V_2\}, \,
    \pi^{-1}(S_4)=\{U_1,U_6,V_1,V_2\}.
  $$
  We will show that the correct one is $\{U_1,U_4,V_1,V_2\}$.

  Let $V_j'=\rho_2(V_j)$ and $V_j''=\rho_2^2(V_j)$ for $j=1,2$. If
  $\pi^{-1}(S_4)=\{U_1,U_2,V_1,V_2\}$, then we have
  $$
    z_4=\frac{Bz_6(1-3z_6)(z_6-z_6(V_1))^4(z_6-z_6(V_2))^4}
    {(1+z_6)(1+3z_6)(z_6-z_6(V_1''))^4(z_6-z_6(V_2''))^4}
  $$
  for some constant $B$. The values of $z_6(V_1)$ and etc. must
  satisfy
  \begin{equation} \label{equation: 444444 Belyi}
  \begin{split}
   &Bx(1-3x)(1-x/z_6(V_1))^4(1-x/z_6(V_2))^4 \\
   &\qquad\qquad\qquad
     -(1+x)(1+3x)(1-x/z_6(V_1''))^4(1-x/z_6(V_2''))^4 \\
   &\qquad=C(1-x)(1-x/z_6(V_1'))^4(1-x/z_6(V_2'))^4
  \end{split}
  \end{equation}
  for some constant $C$. Now if we let
  $p_1(x)=1+ax+bx^2=(1-x/z_6(V_1))(1-x/z_6(V_2))$, then
  $(1-x/z_6(V_1'))(1-x/z_6(V_2'))$ and
  $(1-x/z_6(V_1'')(1-x/z_6(V_2''))$ are scalar multiples of
  \begin{equation*}
  \begin{split}
   p_2(x)&=(1+3x)^2p_1\left(\frac{x-1}{3x+1}\right)=(1-a+b)+(6-2a-2b)x
     +(9+3a+b)x^2, \\
   p_3(x)&=(1-3x)^2p_1\left(\frac{x+1}{1-3x}\right)=(1+a+b)+(-6-2a+2b)z
     +(9-3a+b)x^2,
  \end{split}
  \end{equation*}
  respectively. Substituting these into \eqref{equation: 444444 Belyi}
  and equating the coefficients in the two sides, we find $A=B=0$,
  $a=-2$, $b=-3$, but obviously this is invalid. This means that
  $\pi^{-1}(S_4)\neq\{U_1,U_2,V_1,V_2\}$. Likewise,
  $\pi^{-1}(S_4)\neq\{U_1,U_6,V_1,V_2\}$. Thus, we must have
  $\pi^{-1}(S_4)=\{U_1,U_4,V_1,V_2\}$. Now equating the coefficients
  in the two sides of
  $$
    Bx(1+ax+bx^2)^4-(1-x)(1+3x)p_2(x)^4=C(1+x)(1-3x)p_3(x)^4
  $$
  and excluding the invalid solutions, we get the claimed relation
  between $z_4$ and $z_6$. The relation between $z_3$ and $z_5$ can be
  determined by the known relation between $z_3$ and $z_4$, that between
  $z_4$ and $z_6$, and that between $z_5$ and $z_6$. This process also
  determines the choices of the third roots of unity in the relation
  between $z_3$ and $z_4$ and that between $z_5$ and $z_6$. We omit the
  details.
\end{proof}

\begin{Lemma} \label{lemma: 444444 Schwarzian}
  The automorphic derivative $Q(z_6)=D(z_6,\tau)$ is equal to
  \begin{equation} \label{equation: 444444 Schwarzian}
  \begin{split}
  &\frac{15}{64}\left(\frac1{z_6^2}+\frac1{(1-z_6)^2}+\frac1{(1+z_6)^2}
    +\frac1{(z_6-1/3)^2}+\frac1{(z_6+1/3)^2}\right) \\
  &\qquad\qquad+\frac{45}{128}\left(\frac1{1-z_6}+\frac1{1+z_6}
  +\frac3{1-3z_6}+\frac3{1+3z_6}\right).
  \end{split}
  \end{equation}
\end{Lemma}

\begin{proof} By Proposition \ref{theorem: basis}, the rational function
  $R(x)$ such that automorphic $Q(z_6)=D(z_6,\tau)$ is equal to
  $R(z_6)$ is equal to
  \begin{equation*}
  \begin{split}
  R(x)&=\frac{15}{64}\left(\frac1{x^2}+\frac1{(1-x)^2}+\frac1{(1+x)^2}
    +\frac1{(x-1/3)^2}+\frac1{(x+1/3)^2}\right) \\
  &\qquad\qquad+\frac{B_1}x+\frac{B_2}{x-1}+\frac{B_3}{x+1}+\frac{B_4}{x-1/3}
    +\frac{B_5}{x+1/3}
  \end{split}
  \end{equation*}
  for some constants $B_j$ satisfying
  \begin{equation} \label{equation: 444444 B}
    B_1+B_2+B_3+B_4+B_5=0, \quad
    B_2-B_3+\frac13B_4-\frac13B_5+\frac{15}{16}=0.
  \end{equation}
  Now the normalizer of $\Gamma_6$ in $\SL(2,\R)$ contains at least
  the group of signature $(2,6,8)$. The factor group, in terms of the
  Hauptmodul $z_6$, is generated by
  $\sigma:z_6\mapsto(3z_6+1)/(-3z_6+3)$ and $\tau:z_6\mapsto-z_6$.
  By Proposition \ref{proposition: D 1}, $R(x)$ satisfies
  $$
    R(-x)=R(x), \qquad \frac{144}{(-3x+3)^4}R
    \left(\frac{3x+1}{-3x+3}\right)=R(x).
  $$
  Combining these informations with \eqref{equation: 444444 B}, we
  find
  $$
    B_1=0, \quad B_2=B_4=-\frac{45}{128}, \quad
    B_3=B_5=\frac{45}{128}.
  $$
  This gives us the formula.
\end{proof}

We now prove the theorem.

\begin{proof}[Proof of Theorem \ref{theorem: Class III}]
By Proposition \ref{proposition: dimension}, we have
\begin{equation*}
\begin{split}
  \dim S_6(\Gamma_1)&=1-6+\gauss{\frac62\left(1-\frac14\right)}
    +2\gauss{\frac62\left(1-\frac16\right)}=1, \\
  \dim S_6(\Gamma_4)&=1-6+3\gauss{\frac62\left(1-\frac14\right)}=1, \\
  \dim S_6(\Gamma_6)&=1-6+6\gauss{\frac62\left(1-\frac14\right)}=7.
\end{split}
\end{equation*}
By Proposition \ref{theorem: triangle}, the one-dimensional spaces
$S_6(\Gamma_1)$ and $S_6(\Gamma_4)$ are spanned by
\begin{equation} \label{equation: 444444 F1}
  F_1=z_1^{1/4}(1-z_1)^{1/2}\left({}_2F_1\left(\frac5{24},\frac38;
  \frac34;z_1\right)
  +C_1z_1^{1/4}{}_2F_1\left(\frac{11}{24},\frac58;\frac54;z_1\right)\right)^6
\end{equation}
and
\begin{equation} \label{equation: 444444 F2}
  F_2=z_4^{1/4}(1-z_4)^{1/4}\left({}_2F_1\left(\frac18,\frac38;\frac34;z_4\right)
  +C_2z_4^{1/4}{}_2F_1\left(\frac38,\frac58;\frac54;z_4\right)\right)^6
\end{equation}
for some complex numbers $C_1$ and $C_2$, respectively. Furthermore,
by Proposition \ref{theorem: basis}, if we let
\begin{equation*}
\begin{split}
  f_1&=z_6^{3/8}\left(1-\frac{15}7z_6^2-\frac{111}{14}z_6^4
    -\frac{2045}{46}z_6^6-\frac{11355195}{39928}z_6^8
    -\frac{77997477}{39928}z_6^{10}-\cdots\right) \\
 f_2&=z_6^{5/8}\left(1-\frac53z_6^2-\frac{245}{34}z_6^4
    -\frac{7269}{170}z_6^6-\frac{115223}{408}z_6^8
    -\frac{55230121}{27880}-\cdots\right)
\end{split}
\end{equation*}
be a basis for the solution space of the Schwarzian differential
equation $d^2f/dz_6^2+Q(z_6)f=0$, where $Q(z_6)$ is the rational
function in \eqref{equation: 444444 Schwarzian}, then a basis for
$S_8(\Gamma_6)$ is
$$
  \{z_6^jg:~j=0,\ldots,6\}, \qquad g=\frac{(f_1+C_3f_2)^8}
  {z_6^2(1-z_6^2)^2(1-9z_6^2)^2}.
$$
Now from Lemma \ref{lemma: Class III Belyi}, we have
$$
  z_1=\frac{12\alpha z_6(1-z_6^2)(1-9z_6^2)}{(1+\alpha z_6)^6}
$$
and
$$
  z_4=-\frac{4(1+\beta)^4z_6(1+(-7+4\beta)z_6^2/3)^4}
  {(1+z_6)(1-3z_6)(1+(4+2\beta)z_6-(1+2\beta)z_6^2)^4},
$$
where $\alpha$ is a root of $x^2+3=0$ and $\beta$ is a root of
$x^2+2=0$. Substituting these into \eqref{equation: 444444 F1} and
\eqref{equation: 444444 F2} and comparing the coefficients, we find
$$
  F_1=c_1(1+3z_6^2)^3g
$$
and
\begin{equation*}
\begin{split}
  F_2&=c_2\left(1+\frac{-7+4\beta}3z_6^2\right)
  \left(1+(4+2\beta)z_6-(1+2\beta)z_6^2\right) \\
 &\qquad\qquad\times\left(1-(4+2\beta)z_6-(1+2\beta)z_6^2\right)g
\end{split}
\end{equation*}
for some constants $c_1$ and $c_2$. Taking the sixth roots of $F_1$
and $F_2$ and simplifying, we obtain the identities claimed in the
theorem.
\end{proof}
\end{section}

\begin{section}{Algebraic transformations associated to Class VI}
\label{section: class VI}

According to Appendix A, the subgroup diagram for Takeuchi's Class VI
is
$$
  \begin{diagram}
  \node{(2,4,5)} \arrow{s,l,-}{2} \arrow{ese,r,-}{6}
    \node[4]{(2,4,10)} \arrow{wsw,l,-}{2} \arrow{s,r,-}{2} \\
  \node{(2,5,5)} \arrow{ese,l,-}{6}
    \node[2]{(4,4,5)} \arrow{s,l,-}{2}
    \node[2]{(2,10,10)} \arrow{wsw,l,-}{2} \arrow{s,r,-}{2} \\
  \node[3]{(2,2,5,5)} \arrow{s,l,-}{2}
    \node[2]{(5,10,10)} \arrow{wsw,r,-}{2} \\
  \node[3]{(5,5,5,5)}
  \end{diagram}
$$
Let $\Gamma_1=(2,5,5)$, $\Gamma_2=(5,10,10)$, $\Gamma_3=(5,5,5,5)$,
and $X_1$, $X_2$, $X_3$ be the Shimura curves associated to these
three groups. (The reader is reminded that the subgroup
diagram should be read as ``\emph{there are} arithmetic Fuchsian
subgroups of $\SL(2,\R)$ such that their subgroup relations are given
by the diagram''.) The subgroups relations $\Gamma_3<\Gamma_1,\Gamma_2$
admit Coxeter decompositions as the following figures show.

\centerline{\epsfig{file=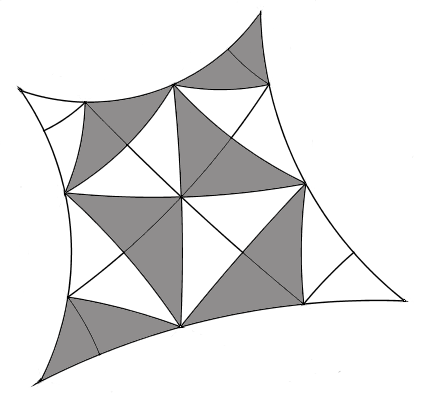,width=1.4in,height=1.3in} \qquad
  \epsfig{file=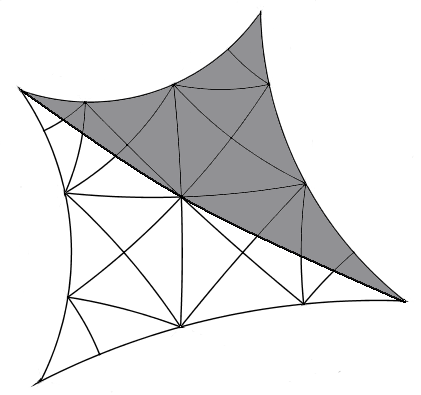,width=1.4in,height=1.3in}}

\noindent Here the small triangles are $(2,4,5)$-triangles. Associated
to this triplet of groups is the following identities.

\begin{Theorem} We have
\begin{equation} \label{equation: favorite 1}
\begin{split}
  &{}_2F_1\left(\frac1{20},\frac14;\frac45;\frac{64z(1-z-z^2)^5}
     {(1-z^2)(1+4z-z^2)^5}\right) \\
  &\qquad\qquad=(1-z^2)^{1/20}(1+4z-z^2)^{1/4}
  {}_2F_1\left(\frac3{10},\frac25;\frac9{10};z^2\right).
\end{split}
\end{equation}
and
\begin{equation} \label{equation: favorite 2}
\begin{split}
  &(1-z-z^2){}_2F_1\left(\frac14,\frac9{20};\frac65;\frac{64z(1-z-z^2)^5}
     {(1-z^2)(1+4z-z^2)^5}\right) \\
  &\qquad\qquad=(1-z^2)^{1/4}(1+4z-z^2)^{5/4}
  {}_2F_1\left(\frac25,\frac12;\frac{11}{10};z^2\right).
\end{split}
\end{equation}
\end{Theorem}

\begin{proof} Label the elliptic points of $X_j$ by $P_2$, $P_5$,
  $P_5'$ for $X_1$, $Q_5$, $Q_{10}$, $Q_{10}'$ for $X_2$, and $R_i$,
  $i=1,\ldots,4$, for $X_3$ such that the ramifications data are given
  by

\centerline{
\begin{tikzpicture}[scale=0.05]
  \node[place] (R1) [label=above:$R_1$, token=1] at (-110,20) {};
  \node[place] (R3) [label=above:$R_3$, token=1] at (-90,20) {};
  \node[place] (S1) [label=above:$S_1$, token=1] at (-70,20) {};
  \node[place] (S3) [label=above:$S_3$, token=1] at (-50,20) {};
  \node[place] (R2) [label=above:$R_2$, token=1] at (-30,20) {};
  \node[place] (R4) [label=above:$R_4$, token=1] at (-10,20) {};
  \node[place] (S2) [label=above:$S_2$, token=1] at (10,20) {};
  \node[place] (S4) [label=above:$S_4$, token=1] at (30,20) {};
  \node[place] (R21) [label=above:$R_2$, token=1] at (50,20) {};
  \node[place] (R41) [label=above:$R_4$, token=1] at (70,20) {};
  \node[place] (R11) [label=above:$R_1$, token=1] at (90,20) {};
  \node[place] (R31) [label=above:$R_3$, token=1] at (110,20) {};
  \node[place] (P51) [label=below:$P_5$, token=1] at (-80,0) {};
  \node[place] (P52) [label=below:$P_5'$, token=1] at (0,0) {};
  \node[place] (Q5) [label=below:$Q_5$, token=1] at (60,0) {};
  \node[place] (Q101) [label=below:$Q_{10}$, token=1] at (90,0) {};
  \node[place] (Q102) [label=below:$Q_{10}'$, token=1] at (110,0) {};
  \draw[line width=.6pt] (R1) --  node[left=2pt] {$1$} (P51);
  \draw[line width=.6pt] (R3) --  node[left=2pt] {$1$}(P51);
  \draw[line width=.6pt] (S1) --  node[right=2pt] {$5$}(P51);
  \draw[line width=.6pt] (S3) --  node[right=2pt] {$5$}(P51);
  \draw[line width=.6pt] (R2) -- node[left=2pt] {$1$}(P52);
  \draw[line width=.6pt] (R4) -- node[left=2pt] {$1$}(P52);
  \draw[line width=.6pt] (S2) -- node[right=2pt] {$5$}(P52);
  \draw[line width=.6pt] (S4) -- node[right=2pt] {$5$}(P52);
  \draw[line width=.6pt] (R21) -- node[left=2pt] {$1$} (Q5);
  \draw[line width=.6pt] (R41) -- node[right=2pt] {$1$} (Q5);
  \draw[line width=.6pt] (R11) -- node[right=2pt] {$2$} (Q101);
  \draw[line width=.6pt] (R31) -- node[right=2pt] {$2$} (Q102);
\end{tikzpicture}
}
\noindent Here the numbers next to the lines are the ramification
indices. We have omitted $P_2$ from the diagram. There are $6$ points
lying above $P_2$. Each has ramficiation index $2$. Choose Hauptmoduls
$z_j$ for $X_j$ by requiring
$$
  z_1(P_5)=0, \ z_1(P_2)=1, \ z_1(P_5')=\infty, \quad
  z_2(Q_{10})=0, \ z_2(Q_5)=1, \ z_2(Q_{10}')=\infty
$$
and
$$
  z_3(R_1)=0, \ z_3(R_2)=1, \ z_3(R_3)=\infty.
$$
The relation between $z_2$ and $z_3$ is easy to figure out.
% The
%ramification data at $Q_{10}$ and $Q_{10}'$ shows that $z_2=Az_3^2$
%for some $A$. Then because $R_2$ is a ramified point, we have $A=4$.
%That is,
We have
\begin{equation} \label{equation: 5555 z2}
  z_2=z_3^2,
\end{equation}
which implies that $z_3(R_4)=-1$. To determine the relation between
$z_1$ and $z_3$, we observe that the tessellation of the hyperbolic
plane by the $(5,5,5,5)$-polygons has extra symmetries by rotation by
$90$ degree around the center of any $(5,5,5,5)$-polygon. In terms of
groups, this means that $\Gamma_3$ has a supergroup $\Gamma$
normalizing $\Gamma_3$ such that $\Gamma/\Gamma_3$ is cyclic of order
$4$. (In fact, $\Gamma$ is the $(4,4,5)$-triangle group in the
subgroup diagram.) Therefore, the automorphism group of $X_3$ has
an element $\sigma$ of order $4$ that permutes $R_1,R_2,R_3,R_4$
cyclically. In terms of the Hauptmodul, we have
$$
  \sigma: z_3\longmapsto\frac{z_3+1}{z_3-1}.
$$
Thus, if the value of $z_3$ at $S_1$ is $a$, then we have
$$
  z_3(S_1)=a, \quad z_3(S_2)=\frac{a-1}{a+1}, \quad
  z_3(S_3)=-\frac{1}{a}, \quad
  z_3(S_4)=-\frac{a+1}{a-1}.
$$
Therefore, the relation between $z_1$ and $z_3$ is
$$
  z_1=\frac{Bz_3(z_3-a)^5(z_3+1/a)^5}
   {(1-z_3^2)(z_3-(a-1)/(a+1))^5(z_3+(a+1)/(a-1))^5}
$$
for some constant $B$. Moreover, the automorphism $\sigma$ of $X_3$
rotates $4$ of the six points lying above $P_2$ cyclically and fixes
the other two. (The reader is reminded that each
$(5,5,5,5)$-polygon represents only half of the fundamental domain for
$\Gamma_3$. The two fixed of $\sigma$ are the centers of the
$(5,5,5,5)$-polygons.) In terms of the Hauptmodul $z_3$, this means
that the values of $z_3$ at the two fixed points of $\sigma$ are $\pm
i$ and if the value of $z_3$ at one of the other $4$ points above
$P_2$ is $b$, then the values at the other $3$ points are
$-1/b$, $(b-1)/(b+1)$, and $-(b+1)/(b-1)$. Thus, we have
\begin{equation*}
\begin{split}
  z_1-1=\frac{C(1+z_3^2)^2(z_3-b)^2(z_3+1/b)^2(z_3-(b-1)/(b+1))^2
    (z_3+(b+1)/(b-1))^2}{(1-z_3^2)(z_3-(a-1)/(a+1))^5(z_3+(a+1)/(a-1))^5}
\end{split}
\end{equation*}
for some constant $C$. Comparing the two sides, we find $a=0,\pm 1,\pm
i$, $a^2+a-1=0$, or $a^2-a-1=0$. The first five solutions are
invalid. The other two solutions give
\begin{equation} \label{equation: 5555 z1 1}
  z_1=\frac{64z_3(1-z_3-z_3^2)^5}{(1-z_3^2)(1+4z_3-z_3^2)^5}
\end{equation}
or
\begin{equation} \label{equation: 5555 z1 2}
  z_1=-\frac{64z_3(1+z_3-z_3^2)^5}{(1-z_3^2)(1-4z_3-z_3^2)^5}.
\end{equation}
Both are valid because of the following reason. Notice that $\Gamma_2$
normalizes $\Gamma_3$. If we take an element $\gamma$ of $\Gamma_2$
not in $\Gamma_3$, then $\gamma^{-1}\Gamma_1\gamma$ is again a
triangle of signature $(2,5,5)$ containing the same $\Gamma_3$.
If the relation between the Hauptmoduls of $\Gamma_1$ and $\Gamma_3$
is \eqref{equation: 5555 z1 1}, then the relation between the
Hauptmoduls of $\gamma^{-1}\Gamma_1\gamma$ and $\Gamma_3$ will be
\eqref{equation: 5555 z1 2}.

Having determining the relations among Hauptmoduls, we can composite
identity (36) in \cite{Vidunas} with Kummer's quadratic transformation
several times to get the identities in the theorem. However, the
procedure is very tedious. Here we provide a better proof using the
theory of automorphic forms on Shimura curves.

By Proposition \ref{proposition: dimension}, we have
\begin{equation*}
\begin{split}
  \dim S_8(\Gamma_1)&=1-8+\gauss{\frac82\left(1-\frac12\right)}
    +2\gauss{\frac82\left(1-\frac15\right)}=1, \\
  \dim S_8(\Gamma_2)&=1-8+\gauss{\frac82\left(1-\frac15\right)}
    +2\gauss{\frac82\left(1-\frac1{10}\right)}=2, \\
  \dim S_8(\Gamma_3)&=1-8+4\gauss{\frac82\left(1-\frac15\right)}=5.
\end{split}
\end{equation*}
By Proposition \ref{theorem: triangle}, the one-dimensional space
$S_8(\Gamma_1)$ is spanned by
\begin{equation} \label{equation: 5555 F1}
  F_1=z_1^{1/5}\left({}_2F_1\left(\frac1{20},\frac14;\frac45;z_1\right)
   +C_1z_1^{1/5}{}_2F_1\left(\frac14,\frac9{20};\frac65;z_1\right)\right)^8
\end{equation}
for some constant $C_1$, and the function
\begin{equation} \label{equation: 5555 F2}
  F_2=z_2^{3/5}(1-z_2)^{1/5}\left(
    {}_2F_1\left(\frac3{10},\frac25;\frac9{10};z_2\right)
  +C_2z_2^{1/10}{}_2F_1\left(\frac25;\frac12;\frac{11}{10};z_2\right)\right)^8
\end{equation}
is contained in $S_8(\Gamma_2)$ for some constant $C_2$. To get a
basis for $S_8(\Gamma_3)$, we need to work out the Schwarzian
differential equation associated to $z_3$. It is actually easy in this
case.

By Proposition \ref{theorem: basis}, the function $z_3'(\tau)$, as a
function of $z_3$, satisfies
$$
  \frac{d^2}{dz_3^2}f+Q(z_3)f=0,
$$
where
\begin{equation} \label{equation: 5555 temp 1}
  Q(z_3)=\frac6{25}\left(\frac1{z_3^2}+\frac1{(1-z_3)^2}+\frac1{(1+z_3)^2}\right)
  +\frac{B_1}{z_3}+\frac{B_2}{z_3-1}+\frac{B_3}{z_3+1}
\end{equation}
for some complex numbers satisfying
\begin{equation} \label{equation: 5555 temp 2}
  B_1+B_2+B_3=0, \quad B_2-B_3+\frac{12}{25}=0.
\end{equation}
To determine the values of $B_j$, we use the automorphism of $X_3$
coming from the normal subgroup relation $\Gamma_3\lhd\Gamma_1$.
Let $\gamma$ be an element of $\Gamma_2$ not in $\Gamma_3$. We know
that
$$
  z_3(\gamma\tau)=-z_3(\tau).
$$
Now by Proposition \ref{proposition: D}, we have
\begin{equation} \label{equation: 5555 temp 3}
\begin{split}
   D(-z_3(\tau),\tau)&=D(z_3(\gamma\tau),\tau)=D(z_3(\gamma\tau),\gamma\tau)
   +D(\gamma\tau,\tau)/(d\gamma\tau/d\tau)^2=Q(z_3(\gamma\tau)) \\
  &=\frac6{25}\left(\frac1{z_3^2}+\frac1{(1-z_3)^2}+\frac1{(1+z_3)^2}\right)
  +\frac{B_1}{-z_3}+\frac{B_2}{-z_3-1}+\frac{B_3}{-z_3+1}
\end{split}
\end{equation}
On the other hand, we also have, by the same proposition,
\begin{equation} \label{equation: 5555 temp 4}
  D(-z_3,\tau)=D(-z_3,z_3)+D(z_3,\tau)/(-1)^2=Q(z_3)
\end{equation}
Comparing \eqref{equation: 5555 temp 1}, \eqref{equation: 5555 temp
  3}, and \eqref{equation: 5555 temp 4}, we find $B_1=0$ and
$B_2=-B_3$. Together with \eqref{equation: 5555 temp 2}, this gives us
$$
  B_1=0, \quad B_2=-\frac6{25}, \quad B_3=\frac6{25}
$$
and
$$
  Q(z_3)=\frac6{25}\left(\frac1{z_3^2}+\frac1{(1-z_3)^2}+\frac1{(1+z_3)^2}
  +\frac1{1-z_3}+\frac1{1+z_3}\right).
$$
Now a basis for the solution space of the Schwarzian differential
equation $d^2f/dz_3^2+Q(z_3)f=0$ is given by
\begin{equation*}
\begin{split}
  f_1&=z_3^{2/5}\left(1-\frac4{15}z_3^2-\frac{52}{475}z_3^4
    -\frac{13436}{206625}z_3^6-\frac{46348}{1033125}z_3^8
    -\frac{2024924}{60265625}z_3^{10}-\cdots\right) \\
  f_2&=z_3^{3/5}\left(1-\frac{12}{55}z_3^2-\frac{28}{275}z_3^4
    -\frac{2708}{42625}z_3^6-\frac{393636}{8738125}z_3^8
    -\frac{7503908}{218453125}z_3^{10}-\cdots\right).
\end{split}
\end{equation*}
By Proposition \ref{theorem: basis},
$$
  g,\ z_3g, \ z_3^2g, \ z_3^3g,\ z_3^4g, \qquad
  g=\frac{(f_1+C_3f_2)^8}{z_3^3(1-z_3)^3(1+z_3)^3},
$$
form a basis for $S_8(\Gamma_3)$ for some constant $C_3$. That is,
after substituting \eqref{equation: 5555 z1 1} and \eqref{equation: 5555 z2}
into \eqref{equation: 5555 F1} and \eqref{equation: 5555 F2},
respectively, we have $F_1=h_1(z_3)g$ and $F_2=h_2(z_3)g$ for some
polynomials $h_1(x)$ and $h_2(x)$ of degree $\le 4$. Indeed, by
comparing the coefficients, we find
$$
  F_1=2^{6/5}(1-z_3-z_3^2)(1+4z_3-z_3^2)g, \qquad
  F_2=z_3g.
$$
(The computation becomes easier if we take the $8$th roots of the
functions first.) Simplifying the relation
$z_3F_1=2^{6/5}(1-z_3-z_3^2)(1+4z_3-z_3^2)F_2$, we get the two
identities in the theorem. This completes the proof.
\end{proof}

\end{section}

\begin{section}{Algebraic transformations associated to other classes}
%  Note that the quaternion algebra in Class I is $M(2,\Q)$, so the
%  Shimura curves are just the classical modular curves. In this case,
%  it is easier to use Fourier expansions of modular forms and modular
%  functions. We will not discuss this case.

\begin{subsection}{Classes II, V, and XII}
The subgroup diagrams of Class II, V, and XII are all of the form
$$
  \begin{diagram}
  \node[3]{(2,4,2n)} \arrow{sw,l,-}{2} \arrow{se,r,-}{2} \\
  \node[2]{(2,2n,2n)} \arrow{sw,l,-}{2} \arrow{se,r,-}{2}
  \node[2]{(4,4,n)} \arrow{sw,r,-}{2} \\
  \node{(n,2n,2n)} \arrow{se,l,-}{2} \node[2]{(2,n,2,n)} \arrow{sw,r,-}{2} \\
  \node[2]{(n,n,n,n)}
  \end{diagram}
$$
The subgroup relation $(2,2n,2n)\cap(4,4,n)=(2,n,2,n)$ is a special
case of
$$
  \begin{diagram}
  \node[2]{(2,2m,2n)} \arrow{sw,l,-}{2} \arrow{se,r,-}{2} \\
  \node{(m,2n,2n)} \arrow{se,l,-}{2} \node[2]{(n,2m,2m)}
  \arrow{sw,r,-}{2} \\
  \node[2]{(m,n,m,n)}
  \end{diagram}
$$
which arises from the Coxeter decompositions of a quadrilateral
polygon that is symmetric with respect to both the diagonals as shown below
$$
  \begin{diagram}
  \node[2]{\epsfig{file=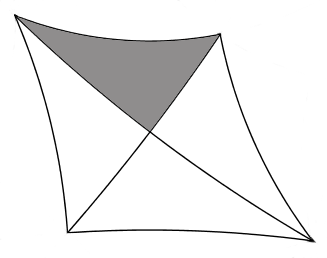,width=0.46in,height=0.36in}}
    \arrow{sw,l,-}{2} \arrow{se,r,-}{2} \\
  \node{\epsfig{file=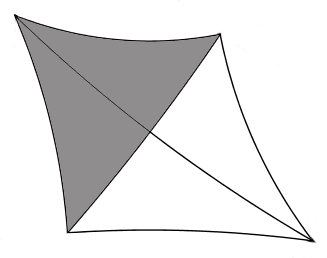,width=0.46in,height=0.36in}}
    \arrow{se,l,-}{2}
  \node[2]{\epsfig{file=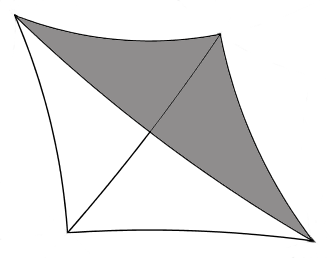,width=0.46in,height=0.36in}}
    \arrow{sw,r,-}{2} \\
  \node[2]{\epsfig{file=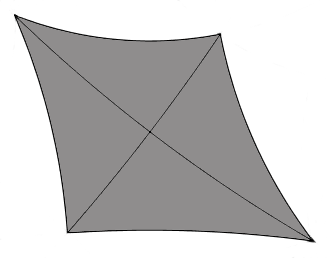,width=0.46in,height=0.36in}}
  \end{diagram}
$$
Associated to this family of subgroup relations is the
following identity.

\begin{Theorem} For real numbers $a$ and $b$ such that neither $b+3/4$
  nor $2b+1/2$ is a nonpositive integer, we have
  \begin{equation*}
    (1+z)^{2a+2b}{}_2F_1\left(a+b,a+\frac14;b+\frac34;z^2\right)
   ={}_2F_1\left(a+b,b+\frac14;2b+\frac12;\frac{4z}{(1+z)^2}\right)
  \end{equation*}
  in a neighborhood of $z=0$.
\end{Theorem}

This identity can be easily proved using Kummer's quadratic
transformation. Alternatively, one can verify that both sides are
solutions of the differential equation
$$
  2z(1-z)(1+z)^2F''-(1+z)((3-4b)z^2+8(a+b)z-4b-1)F'
  -(a+b)(1+4b)(1-z)F=0.
$$
and that the local behaviors at $z=0$ agree. We omit the details.
\end{subsection}

\begin{subsection}{Classes IV, VIII, XI, XIII, XV, XVII}

The subgroups diagrams of Classes IV, VIII, XI, XIII, XV, and XVIII
are either of the form
$$ \tiny
  \begin{diagram}
  \node[3]{(2,3,12n)} \arrow{wsw,l,-}{2} \arrow{s,l,-}{4} \arrow{ese,r,-}{3} \\
  \node{(3,3,6n)} \arrow{s,l,-}{4} \arrow{ese,r,-}{} \node[2]{(3,4n,12n)}
    \arrow{wsw,r,-}{} \arrow{ese,r,-}{} \node[2]{(2,6n,12n)}
    \arrow{wsw,r,-}{} \arrow{s,r,-}{4} \arrow{ese,r,-}{2} \\
  \node{(3,3,2n,6n)} \arrow{ese,l,-}{3} \node[2]{(6n,6n,6n)}
    \arrow{s,l,-}{4} \arrow{ese,r,-}{2} \node[2]{(2n,4n,6n,12n)}
    \arrow{wsw,l,-}{2} \arrow{ese,r,-}{} \node[2]{(3n,12n,12n)}
    \arrow{wsw,r,-}{2} \arrow{s,r,-}{4} \\
  \node[3]{(2n^3,6n^3)} \arrow{ese,l,-}{2} \node[2]{(3n^2,6n^2)}
    \arrow{s,l,-}{4} \node[2]{(n,3n,4n^2,12n^2)} \arrow{wsw,r,-}{2} \\
  \node[5]{(1;n^2,2n^2,3n^2,6n^2)}
  \end{diagram}
$$
or sub-diagram of it with Class XI having one extra node. There are two
families of essentially new identities associated to these classes.
One corresponds to the pair of $(3,3,6n)$ and $(3,4n,12n)$. (Theorem
\ref{theorem: 3-3-6n} below.) One corresponds to the pair of
$(3,4n,12n)$ and $(2,6n,12n)$. (Theorem \ref{theorem: 2-6n-12n}
below.)
%All the other identities associated to these classes are
%either classical results or the composition of classical results and
%the identities in this section and the previous section.

\begin{Theorem} \label{theorem: 3-3-6n}
  For a real number $a$ such that neither $3a+1$ nor $2a+1$ is a
  nonpositive integer, we have
  \begin{equation*}
  \begin{split}
   &(1+z)^{a+1/6}(1-z/3)^{3a+1/2}{}_2F_1\left(2a+\frac13,
    a+\frac13;3a+1;z^2\right) \\
  &\qquad\qquad={}_2F_1\left(a+\frac16,a+\frac12;2a+1;
   \frac{16z^3}{(1+z)(3-z)^3}\right)
  \end{split}
  \end{equation*}
  in a neighborhood of $z=0$.
\end{Theorem}

\begin{Theorem} \label{theorem: 2-6n-12n}
  For a real number $a$ such that neither $6a+1$ nor $4a+1$ is a
  nonpositive integer, we have
  \begin{equation*}
  \begin{split}
   &(1-z)^{9a+3/4}{}_2F_1\left(4a+\frac13,2a+\frac13;6a+1;
    -\frac{27z^2(1-z)}{1-9z}\right) \\
   &\qquad\qquad=(1-9z)^{a+1/12}{}_2F_1\left(3a+\frac14,a+\frac14;4a+1;
    -\frac{64z^3}{(1-z)^3(1-9z)}\right)
  \end{split}
  \end{equation*}
  in a neighborhood of $z=0$.
\end{Theorem}

In principle, these two identities can be deduced from Kummer's and
Goursat's transformations, once the related Belyi functions are
determined. Here we briefly indicate how one can prove the theorems
in the cases where the parameters correspond to discrete Fuchsian
groups using theory of automorphic forms.

\begin{proof}[Proof of Theorem \ref{theorem: 3-3-6n} in the cases of
  Shimura curves]

For the pair of $(3,3,6n)$ and $(3,4n,12n)$, the subgroup relations
admit Coxeter decompositions, as shown in the figures

\centerline{\epsfig{file=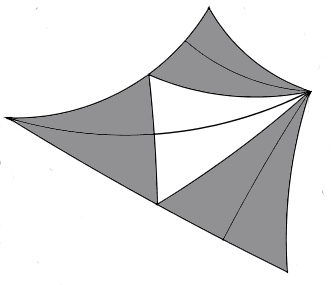,height=1.22in,width=1.44in} \qquad\qquad
\epsfig{file=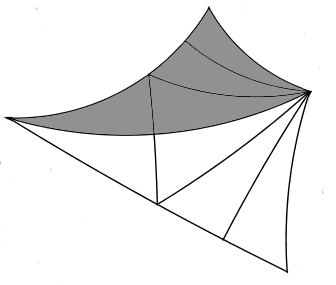,height=1.22in,width=1.44in}}

\noindent Here the parameter $n$ in the figures is $1$ and the smaller
triangles are $(2,3,12)$-triangles. Let $\Gamma_1=(3,3,6n)$,
$\Gamma_2=(3,4n,12n)$, $\Gamma_3=\Gamma_1\cap\Gamma_2$, and let $X_i$,
$i=1,\ldots,3$ be the associated Shimura curves. Denote by $P_3$,
$P_3'$, and $P_{6n}$ the elliptic points of orders $3$, $3$, and $6n$
on $X_1$, by $Q_3$, $Q_{4n}$, and $Q_{12n}$ the elliptic points of
orders $3$, $4n$, and $12n$ on $X_2$, and by $R_3$, $R_3'$, $R_{2n}$,
and $R_{6n}$ the elliptic points of order $3$, $3$, $2n$, and $6n$ on
$X_3$. The points are labelled in a way such that the ramification
data are given by
\centerline{
\begin{tikzpicture}[scale=0.06]
  \node[place] (R31) [label=above:$R_3$, token=1] at (-90,20) {};
  \node[place] (S1) [label=above:$S_1$, token=1] at (-70,20) {};
  \node[place] (R32) [label=above:$R_3'$, token=1] at (-50,20) {};
  \node[place] (S2) [label=above:$S_2$, token=1] at (-30,20) {};
  \node[place] (R2n) [label=above:$R_{2n}$, token=1] at (-10,20) {};
  \node[place] (R6n) [label=above:$R_{6n}$, token=1] at (10,20) {};
  \node[place] (R311) [label=above:$R_3$, token=1] at (30,20) {};
  \node[place] (R321) [label=above:$R_3'$, token=1] at (50,20) {};
  \node[place] (R2n1) [label=above:$R_{2n}$, token=1] at (70,20) {};
  \node[place] (R6n1) [label=above:$R_{6n}$, token=1] at (90,20) {};
  \node[place] (P31) [label=below:$P_3$, token=1] at (-80,0) {};
  \node[place] (P32) [label=below:$P_3'$, token=1] at (-40,0) {};
  \node[place] (P6n) [label=below:$P_{6n}$, token=1] at (0,0) {};
  \node[place] (Q3) [label=below:$Q_3$, token=1] at (40,0) {};
  \node[place] (Q4n) [label=below:$Q_{4n}$, token=1] at (70,0) {};
  \node[place] (Q12n) [label=below:$Q_{12n}$, token=1] at (90,0) {};
  \draw[line width=.6pt] (R31) --  node[left=2pt] {$1$} (P31);
  \draw[line width=.6pt] (S1) --  node[right=2pt] {$3$}(P31);
  \draw[line width=.6pt] (R32) --  node[left=2pt] {$1$}(P32);
  \draw[line width=.6pt] (S2) --  node[right=2pt] {$3$}(P32);
  \draw[line width=.6pt] (R2n) -- node[left=2pt] {$3$}(P6n);
  \draw[line width=.6pt] (R6n) -- node[right=2pt] {$1$}(P6n);
  \draw[line width=.6pt] (R311) -- node[left=2pt] {$1$}(Q3);
  \draw[line width=.6pt] (R321) -- node[right=2pt] {$1$}(Q3);
  \draw[line width=.6pt] (R2n1) -- node[left=2pt] {$2$} (Q4n);
  \draw[line width=.6pt] (R6n1) -- node[right=2pt] {$2$} (Q12n);
\end{tikzpicture}
}
Choose Hauptmoduls $z_j$ on $X_j$, $j=1,2,3$, by requiring
$$
  z_1(P_{6n})=0, \ z_1(P_3)=1, \ z_1(P_3')=\infty, \quad
  z_2(Q_{4n})=0, \ z_2(Q_3)=1, \ z_2(Q_{12n})=\infty
$$
and
$$
  z_3(R_{2n})=0, \ z_3(R_3)=1, \ z_3(R_{6n})=\infty.
$$
It is easy to see from the ramification informations that
\begin{equation} \label{equation: 336 z2}
  z_2=z_3^2,
\end{equation}
which implies that $z_3(R_3')=-1$. For $z_1$, we have
$$
  z_1=\frac{Az_3^3}{(1+z_3)(1-az_3)^3}
$$
for some complex numbers $A$ and $a$, where $1/a$ is the value of
$z_3$ at $S_1$. These two numbers satisfy
\begin{equation} \label{equation: 336 temp}
  1-\frac{Az_3^3}{(1+z_3)(1-az_3)^3}=1-z_1
   =\frac{(1-z_3)(1-bz_3)^3}{(1+z_3)(1-az_3)^3},
\end{equation}
where $1/b$ is the value of $z_3$ at $S_2$. Now observe that
$\Gamma_3$ is a normal subgroup of $\Gamma_2$. Thus, an element of
$\Gamma_2$ not in $\Gamma_3$ induces an automorphism on $X_3$. In
terms of the Hauptmodul $z_3$, it is easy to see that this
automorphism sends $z_3$ to $-z_3$. Since this automorphism maps $S_1$
to $S_2$, we find $b=-a$. Then comparing the two sides of
\eqref{equation: 336 temp}, we get $A=16/27$, $a=1/3$, and
\begin{equation} \label{equation: 336 z1}
  z_1=\frac{16z_3^3}{(1+z_3)(3-z_3)^3}.
\end{equation}
Now by Proposition \ref{proposition: dimension}, we have
\begin{equation*}
\begin{split}
  \dim S_6(\Gamma_1)=\dim S_6(\Gamma_2)=1, \quad
  \dim S_6(\Gamma_3)=\begin{cases}
  2, &\text{if }n=1, \\
  3, &\text{if }n\ge 2. \end{cases}
\end{split}
\end{equation*}
From now on, we assume that $n\ge 2$.

By Proposition \ref{theorem: triangle}, the one-dimensional spaces
$S_6(\Gamma_1)$ and $S_6(\Gamma_2)$ are spanned by
\begin{equation} \label{equation: 336 F1}
\begin{split}
  F_1&=z_1^{1-1/2n}\Bigg({}_2F_1\left(\frac16-\frac1{12n},
    \frac12-\frac1{12n};1-\frac1{6n};z_1\right) \\
  &\qquad\qquad+C_1z_1^{1/6n}{}_2F_1\left(\frac16+\frac1{12n},
    \frac12+\frac1{12n};1+\frac1{6n};z_1\right)\Bigg)^6
\end{split}
\end{equation}
and
\begin{equation} \label{equation: 336 F2}
\begin{split}
  F_2&=z_2^{1-3/4n}\Bigg({}_2F_1\left(\frac13-\frac1{6n},
    \frac13-\frac1{12n};1-\frac1{4n};z_2\right) \\
  &\qquad\qquad+C_2z_2^{1/4n}{}_2F_1\left(\frac13+\frac1{12n},
   \frac13+\frac1{6n};1+\frac1{4n};z_2\right)\Bigg)^6,
\end{split}
\end{equation}
respectively, for some constants $C_1$ and $C_2$. Also, if we let $f_1=z_3^{1/2-1/4n}(1+c_1z+\cdots)$ and
$f_2=z_3^{1/2+1/4n}(1+d_1z+\cdots)$ be a basis of the solution space
of the Schwarzian differential equation $d^2f/dz_3^2+Q(z_3)f=0$
associated to $z_3$, then by Proposition \ref{theorem: basis},
$S_6(\Gamma_3)$ is spanned by $g$, $z_3g$, and $z_3^2g$, where
$$
  g=\frac{(f_1+C_3f_2)^6}{z_3^2(1-z_3)^2(1+z_3)^2}
$$
for some constant $C_3$. Now we substitute \eqref{equation: 336 z1}
and \eqref{equation: 336 z2} into \eqref{equation: 336 F1} and
\eqref{equation: 336 z2}, respectively. We find
$$
  F_1=a_1z_3^{3-3/2n}+\cdots, \qquad F_2=z_3^{2-3/2n}+\cdots,
$$
where $a_1=(16/27)^{1-1/2n}$, and thus
$$
  F_1=a_1z_3^2g, \qquad F_2=(z_3+a_2z_3^2)g
$$
for some constant $a_2$. That is, $a_1zF_2/F_1=1+a_2z_3$. We then take
the $6$th roots of the two sides and compare the coefficients of
$z^{3/2-1/4n}$, we find that $a_2$ is actually $0$. After simplifying,
we arrive at
\begin{equation*}
\begin{split}
  &(1+z)^{1/6-1/12n}(1-z/3)^{1/2-1/4n}
   {}_2F_1\left(\frac13-\frac1{6n},\frac13-\frac1{12n};
   1-\frac1{4n};z^2\right) \\
  &\qquad\qquad={}_2F_1\left(\frac16-\frac1{12n},\frac12-\frac1{12n};
   1-\frac1{6n};\frac{16z^3}{(1+z)(3-z)^3}\right).
\end{split}
\end{equation*}
This proves Theorem \ref{theorem: 3-3-6n} in the case the parameters
correspond to arithmetic triangle groups.
\end{proof}

\begin{proof}[Proof of Theorem \ref{theorem: 2-6n-12n} in the cases of
  Shimura curves]

The subgroups $(3,4n,12n)$, $(2,6n,12n)$ and their intersection admit
Coxeter decompositions as the figures below show.

\centerline{\epsfig{file=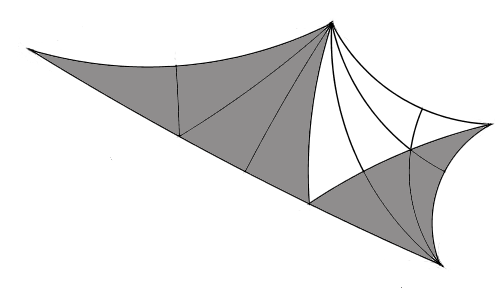,height=1.152in,width=1.976in}
  \qquad \epsfig{file=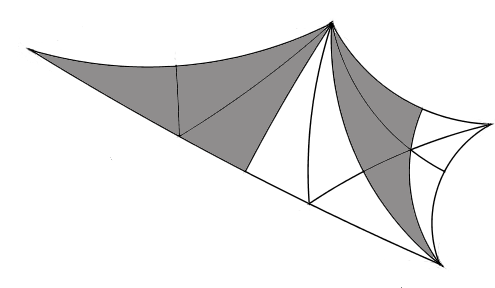,height=1.152in,width=1.976in} }

\noindent Here the parameter $n$ in the figures is $1$ and the small
triangles are $(2,3,12)$-triangles.

Denote the groups $(3,4n,12n)$, $(2,6n,12n)$, and $(2n,4n,6n,12n)$ by
$\Gamma_1$, $\Gamma_2$, and $\Gamma_3$, respectively.
Label the elliptic points of $(3,4n,12n)$ by $P_3$, $P_{4n}$, and
$P_{12n}$, those of $(2,6n,12n)$ by $Q_2$, $Q_{6n}$, and $Q_{12n}$,
and those of $(2n,4n,6n,12n)$ by $R_{2n}$, $R_{4n}$, $R_{6n}$, and
$R_{12n}$. The ramifications are shown as follows.

\centerline{
\begin{tikzpicture}[scale=0.05]
  \node[place] (R2n1) [label=above:$R_{2n}$, token=1] at (-105,20) {};
  \node[place] (R4n1) [label=above:$R_{4n}$, token=1] at (-85,20) {};
  \node[place] (R6n1) [label=above:$R_{6n}$, token=1] at (-65,20) {};
  \node[place] (R12n1) [label=above:$R_{12n}$, token=1] at (-45,20) {};
  \node[place] (P31) [token=1] at (-25,20) {};
  \node[place] (R2n2) [label=above:$R_{2n}$, token=1] at (5,20) {};
  \node[place] (R6n2) [label=above:$R_{6n}$, token=1] at (25,20) {};
  \node[place] (R4n2) [label=above:$R_{4n}$, token=1] at (45,20) {};
  \node[place] (R12n2) [label=above:$R_{12n}$, token=1] at (65,20) {};
  \node[place] (Q21) [token=1] at (85,20) {};
  \node[place] (Q22) [token=1] at (105,20) {};
  \node[place] (P4n) [label=below:$P_{4n}$, token=1] at (-95,0) {};
  \node[place] (P12n) [label=below:$P_{12n}$, token=1] at (-55,0) {};
  \node[place] (P3) [label=below:$P_3$, token=1] at (-25,0) {};
  \node[place] (Q6n) [label=below:$Q_{6n}$, token=1] at (15,0) {};
  \node[place] (Q12n) [label=below:$Q_{12n}$, token=1] at (55,0) {};
  \node[place] (Q2) [label=below:$Q_2$, token=1] at (95,0) {};
  \draw[line width=.6pt] (R2n1) --  node[left=2pt] {$2$} (P4n);
  \draw[line width=.6pt] (R4n1) --  node[right=2pt] {$1$} (P4n);
  \draw[line width=.6pt] (R6n1) --  node[left=2pt] {$2$} (P12n);
  \draw[line width=.6pt] (R12n1) --  node[right=2pt] {$1$}(P12n);
  \draw[line width=.6pt] (P31) -- node[left=2pt] {$3$}(P3);
  \draw[line width=.6pt] (R2n2) -- node[left=2pt] {$3$}(Q6n);
  \draw[line width=.6pt] (R6n2) -- node[right=2pt] {$1$}(Q6n);
  \draw[line width=.6pt] (R4n2) -- node[left=2pt] {$3$}(Q12n);
  \draw[line width=.6pt] (R12n2) -- node[right=2pt] {$1$} (Q12n);
  \draw[line width=.6pt] (Q21) -- node[left=2pt] {$2$} (Q2);
  \draw[line width=.6pt] (Q22) -- node[right=2pt] {$2$} (Q2);
\end{tikzpicture}
}

\noindent Choose Hauptmoduls $z_j$ for $\Gamma_j$, $j=1,\ldots,3$, by
requiring that
\begin{equation*}
\begin{split}
  z_1(P_{4n})=0, \quad &z_1(P_3)=1, \quad z_1(P_{12n})=\infty, \\
  z_2(Q_{6n})=0, \quad &z_2(Q_2)=1, \quad z_2(Q_{12n})=\infty, \\
  z_3(R_{2n})=0, \quad &z_3(R_{4n})=1, \quad z_3(R_{6n})=\infty.
\end{split}
\end{equation*}
It is easy to work out the relation between $z_1$ and $z_3$ and that
between $z_2$ and $z_3$. They are
\begin{equation} \label{equation: 2-6n-12n Belyi}
  z_1=\frac{27z_3^2(1-z_3)}{1-9z_3}, \qquad
  z_2=-\frac{64z_3^3}{(1-z_3)^3(1-9z_3)}.
\end{equation}
Here $1/9$ is the value of $z_3$ at $R_{12n}$. We then follow the same
arguments as before to obtain the claimed identities. We omit the details.
\end{proof}

%By Proposition \ref{proposition: dimension}, we have
%$$
%  \dim S_8(\Gamma_1)=\dim S_8(\Gamma_2)=1, \quad
%  \dim S_8(\Gamma_3)=\begin{cases}
%  4, &\text{if }n=1, \\
%  5, &\text{if }n\ge 2. \end{cases}
%$$
%From now on, we assume that $n\ge 2$.

%By Theorem \ref{theorem: triangle}, the spaces $S_8(\Gamma_1)$ and
%$S_8(\Gamma_2)$ are spanned by
%\begin{equation} \label{equation: 2-6n-12n F1}
%\begin{split}
%  F_1&=z_1^{1-1/n}(1-z_1)^{2/3}\Bigg(
%     {}_2F_1\left(\frac13-\frac1{6n},\frac13-\frac1{12n};
%     1-\frac1{4n};z_1\right) \\
%  &\qquad\qquad\qquad+C_1z_1^{1/4n}
%     {}_2F_1\left(\frac13+\frac1{12n},\frac13+\frac1{6n};
%     1+\frac1{4n};z_1\right)\Bigg)^8
%\end{split}
%\end{equation}
%and
%\begin{equation} \label{equation: 2-6n-12n F2}
%\begin{split}
%  F_2&=z_2^{1-2/3n}\Bigg(
%    {}_2F_1\left(\frac14-\frac1{8n},\frac14-\frac1{24n};
%    1-\frac1{6n};z_2\right) \\
%  &\qquad\qquad\qquad+C_2z_2^{1/6n}
%    {}_2F_1\left(\frac14+\frac1{24n},\frac14+\frac1{8n};
%    1+\frac1{6n};z_2\right)\Bigg)^8.
%\end{split}
%\end{equation}
%\end{proof}
\end{subsection}
\end{section}
\bigskip

\centerline{{\sc Appendix A. List of arithmetic triangle groups}}
\medskip

In this section, we determine the signatures of the intersections of
commensurable triangle groups.

According to \cite{Takeuchi1,Takeuchi2}, there are totally $85$
arithmetic triangle groups, falling in $19$ different commensurability
classes. Here we give the subgroup diagrams. Note that since most
groups here have genus $0$, we omit the genus information from the
signature, unless the group has a positive genus. Also, to save space,
the notation $(g;e_1^{n_1},\ldots,e_r^{n_r})$ means that the Shimura
curve has $n_i$ elliptic points order $e_i$. Furthermore, for
convenience, we will often call the groups by their signatures, even
though this raises some ambiguity.

\begin{Remark} There is some ambiguity when we say ``the
  intersections of commensurable triangle groups'' because there may
  be more than one orders whose norm-one groups have the same
  signature and the intersections of these groups with another group
  may have different signatures. For example, in the case $B=M(2,\Q)$,
  the subgroups $\Gamma_0(2)$ and $\Gamma^0(2)$ of $\SL(2,\Z)$ have
  the same signature $(0;2,\infty,\infty)$ and the group $\Gamma_0(4)$
  has signature $(0;\infty,\infty,\infty)$. The intersection of
  $\Gamma_0(2)$ and $\Gamma_0(4)$ is just $\Gamma_0(4)$, but the
  intersection of $\Gamma^0(2)$ and $\Gamma_0(4)$ has signature
  $(0;\infty,\infty,\infty,\infty)$. Thus, the subgroup diagrams
  described here should be read as ``\emph{there
    are} arithmetic groups whose subgroup relations are given by the
  subgroup diagrams''.
\end{Remark}

Since it is not easy to describe explicitly the orders associated to
arithmetic triangle groups, here we use group theory and properties of
discrete subgroups of $\SL(2,\R)$ to determine the signatures.
We will work out the case of Class IV in \cite{Takeuchi2} and omit the
proof of the others.

According to \cite{Takeuchi2}, Class IV of arithmetic triangle groups
has the following subgroup diagram.
$$
  \begin{diagram}
  \node[3]{(2,3,12)} \arrow{wsw,l,-}{2} \arrow{s,l,-}{4} \arrow{ese,r,-}{3} \\
  \node{(3,3,6)} \arrow{ese,r,-}{3} \node[2]{(3,4,12)}
    \node[2]{(2,6,12)}
    \arrow{wsw,r,-}{2} \arrow{ese,r,-}{2} \\
  \node[3]{(6,6,6)} \node[4]{(3,12,12)}
 \end{diagram}
$$
Here the numbers next to the lines are the indices. Set
\begin{equation*}
\begin{split}
  \Gamma_1&=(2,3,12), \quad \Gamma_2=(3,3,6), \quad \Gamma_3=(3,4,12),
  \\
  \Gamma_4&=(2,6,12), \quad \Gamma_5=(6,6,6), \quad \Gamma_6=(3,12,12),
\end{split}
\end{equation*}
and let $X_i$, $i=1,\ldots,6$, denote the respective Shimura curves.
To determine $\Gamma_2\cap\Gamma_3$, we observe that $\Gamma_2$ is a
normal subgroup of $\Gamma_1$ of index $2$ and
$\Gamma_1=\Gamma_2\Gamma_3$. Thus, $\Gamma_2\cap\Gamma_3$ is a normal
subgroup of $\Gamma_3$ of index $2$. Now the elliptic point of order
$3$ on $X_3$ must split into two points in $X(\Gamma_2\cap\Gamma_3)$
because $2\nmid 3$. Then from the Riemann-Hurwitz formula, we see that
the elliptic points of order $4$ and $12$ must be ramified. That is,
the curve $X(\Gamma_2\cap\Gamma_3)$ must have signature $(2,3,3,6)$.
In fact, this can also be seen from the following figures.

\centerline{\epsfig{file=336.png,height=1.22in,width=1.44in} \qquad\qquad
\epsfig{file=3-4-12-2.png,height=1.22in,width=1.44in}}

\noindent Here the smaller triangles are $(2,3,12)$-triangles. The
figures show that the triangle group $(2,3,12)$ contains two subgroups
of signatures $(3,3,6)$ and $(3,4,12)$, respectively, whose
intersection has signature $(2,3,3,6)$. (In fact, the theoretical
argument above shows that for any pair of subgroups of $\Gamma_1$ with
signatures $(3,3,6)$ and $(3,4,12)$, respectively, the intersection
must have signature $(2,3,3,6)$.)

Likewise, the figures

\centerline{\epsfig{file=2-6-12.png,height=1.152in,width=1.976in} \qquad
\epsfig{file=3-4-12.png,height=1.152in,width=1.976in}}

\noindent show that there are two subgroups of $\Gamma_1$ of signatures
$(2,6,12)$ and $(3,4,12)$ such that there intersection has signature
$(2,4,6,12)$. We have the following subgroup diagram.
$$
  \begin{diagram}
  \node[3]{(2,3,12)} \arrow{wsw,l,-}{2} \arrow{s,l,-}{4} \arrow{ese,r,-}{3} \\
  \node{(3,3,6)} \arrow{s,l,-}{4} \arrow{ese,r,-}{} \node[2]{(3,4,12)}
    \arrow{wsw,-} \arrow{ese,-} \node[2]{(2,6,12)}
    \arrow{wsw,r,-}{} \arrow{ese,r,-}{2} \arrow{s,r,-}{4} \\
  \node{(2,3,3,6)}
  \node[2]{(6,6,6)} \node[2]{(2,4,6,12)} \node[2]{(3,12,12)}
 \end{diagram}
$$
Let $\Gamma_7=(2,3,3,6)$ and $\Gamma_8=(2,4,6,12)$ and $X_7$ and $X_8$
be their associated Shimura curves. Again, because
$\Gamma_5$ is a normal subgroup of $\Gamma_4$ of index $2$ and
$\Gamma_5\Gamma_8=\Gamma_4$, the intersection of $\Gamma_5$ and
$\Gamma_8$ is a subgroup of index $2$ of $\Gamma_8$. Now the group
$(2,4,6,12)$ has many subgroups of index $2$. (The structure of the
quotient group of $(2,4,6,12)$ over its commutator subgroup is
$C_2\times C_4\times C_6$.) To determine which of them is contained is
the group $(6,6,6)$, we use the following properties.
\begin{enumerate}
\item If $p$ is an elliptic point of order $e$ on $X_8$, then its
  preimage in the covering $X(\Gamma_5\cap\Gamma_8)\to X_8$ consists
  of either a single elliptic point of order $e/2$ or two elliptic
  points of order $e$.
\item The total branch number of any finite covering of compact
  Riemann surface is always even.
\item The volume of $X(\Gamma_5\cap\Gamma_8)$ is twice of that of
  $X_8$. Thus, if $(g;e_1,\ldots,e_r)$ is the signature of
  $X(\Gamma_5\cap\Gamma_8)$, then we must have
  $$
    2g-2+\sum_{i=1}^r\left(1-\frac1{e_j}\right)
   =2\left(2-\frac12-\frac14-\frac16-\frac1{12}\right)=2.
  $$
\end{enumerate}
From these informations, we find that possible signatures of a subgroup
of index $2$ of $(2,4,6,12)$ are
\begin{equation} \label{equation: 2-3-12 1}
\begin{split}
 &(1;2,3,6),\ (0;2,6^2,12^2), \ (0;3,4^2,12^2), \
  (0;4^2,6^3) \\
 &\quad(0;2^3,3,12^2), \ (0;2^3,6^3), \ (0;2^2,3,4^2,6).
\end{split}
\end{equation}
Likewise, an elliptic point of order $6$ on $X_5$ can
\begin{enumerate}
\item splits into $4$ elliptic points of order $6$, or
\item splits into $2$ elliptic points of order $3$, or
\item splits into $1$ elliptic point of order $3$ and $2$ elliptic
  point of order $6$, or
\item splits into $1$ elliptic point of order $2$ and $1$ elliptic
  point of order $6$,
\end{enumerate}
in the covering $X(\Gamma_5\cap\Gamma_8)\to X_5$ of degree $4$.
Also, the total branch number of $X(\Gamma_5\cap\Gamma_8)\to X_5$ must
be a positive even integer and the volume of $X(\Gamma\cap\Gamma_8)$
is $2$. We find the possible signatures of a subgroup of index $4$ of
$\Gamma_5$ are
\begin{equation} \label{equation: 2-3-12 2}
  (0;2^3,6^3), \ (0;2^2,3^2,6^2), \ (0;2,3^4,6), \
  (0;3^6).
\end{equation}
From \eqref{equation: 2-3-12 1} and \eqref{equation: 2-3-12 2}, we
conclude that the signature of $\Gamma_5\cap\Gamma_8$ must be
$(0;2^3,6^3)$. This can also be seen from the figures.

\centerline{\epsfig{file=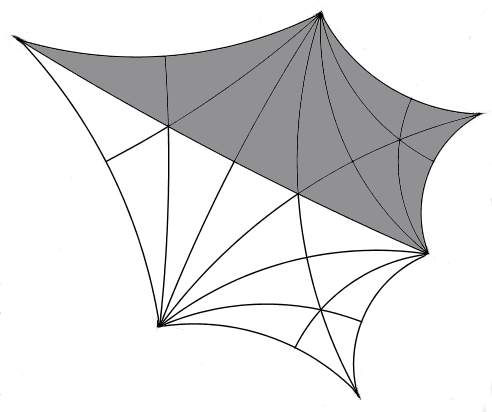,width=1.44in,height=1.2in}
\qquad\qquad \epsfig{file=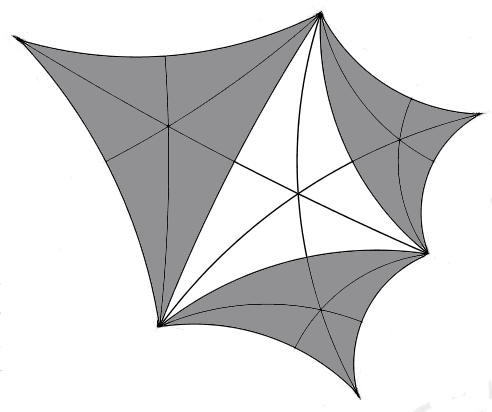,width=1.44in,height=1.2in}}

By the same argument, we can also show that the intersection of
$\Gamma_6$ and $\Gamma_8$ must have signature $(0;3,4^2,12^2)$ and
the intersection of $\Gamma_5$ and $\Gamma_6$ has signature
$(0;3,3,6,6)$. The subgroup diagram becomes
$$
  \begin{diagram}
  \node[3]{(2,3,12)} \arrow{wsw,l,-}{2} \arrow{s,l,-}{4} \arrow{ese,r,-}{3} \\
  \node{(3,3,6)} \arrow{s,l,-}{4} \arrow{ese,r,-}{} \node[2]{(3,4,12)}
    \arrow{wsw,-} \arrow{ese,-} \node[2]{(2,6,12)}
    \arrow{wsw,r,-}{} \arrow{ese,r,-}{2} \arrow{s,r,-}{4} \\
  \node{(2,3,3,6)} \arrow{ese,r,-}{3}
  \node[2]{(6,6,6)} \arrow{s,l,-}{4} \arrow{ese,l,-}{2} 
  \node[2]{(2,4,6,12)} \arrow{wsw,r,-}{} \arrow{ese,l,-}{}
  \node[2]{(3,12,12)} \arrow{s,r,-}{4} \arrow{wsw,r,-}{} \\
  \node[3]{(2^3,6^3)} \node[2]{(3^2,6^2)} \node[2]{(3,4^2,12^2)}
 \end{diagram}
$$
Finally, we can show that the only possible signatures of subgroups of
index $2$ in $(2^3,6^3)$ are
$$
  (0;2^6,3^2,6^2), \ (0;2^4,3,6^4), \ (0;2^2,6^6), \
  (1;2^4,3^3), \ (1;2^2,3^2,6^2), \ (1;3,6^4), \ (2;3^3),
$$
while the only possible signatures of subgroups of index $2$ in
$(3,4^2,12^2)$ are
$$
  (0;3^2,4^4,6^2), \ (0;2,3^2,4^2,6,12^2), \ (0;2^2,3^2,12^4),
  (1;2^2,3^2,6^2).
$$
From these, we see that the common intersection of $(2^3,6^3)$,
$(3,4^2,12^2)$, and $(3^2,6^2)$ has signature $(1;2^2,3^2,6^2)$. This
completes the proof of the case of Class IV.

%\begin{Remark} In literature \cite{}, the decompositions of hyperbolic
%  polygons shown in the figures above are called \emph{Coxeter
%    decompositions}. In general, a Coxeter decomposition is a
%  decomposition of a polygon into finitely many Coxeter polygons such
%  that if two Coxeter polygons share a common side, then they are
%  symmetric with respect to the common side.

%  Note that not all subgroup relations given in Appendix A admit
%  Coxeter decomposition. For example, in Class III, the group
%  $(2,4,8)$ is a subgroup of index $3$ of the group $(2,3,8)$, but
%  there is no way one can decompose a $(2,4,8)$-triangle into a union
%  of three $(2,3,8)$-triangles. In the case of Class IV discussed
%  above, the subgroup relation $(2^3,6^3)<(2,3,3,6)$ does not admit a
%  Coxeter decomposition either.
%\end{Remark}

Now we give the subgroup diagrams for arithmetic triangle groups.
\medskip

\paragraph{\bf Class II}
% $K=\Q$ with $D=6$,
% $\Gamma^{(1)}=\Gamma^{(+)}=(0;2,3,2,3)$, and $\Gamma^{(*)}=(2,4,6)$.
$$
  \begin{diagram}
  \node[3]{(2,4,6)} \arrow{sw,l,-}{2} \arrow{se,r,-}{2} \\
  \node[2]{(2,6,6)} \arrow{sw,l,-}{2} \arrow{se,r,-}{2}
  \node[2]{(3,4,4)} \arrow{sw,r,-}{2} \\
  \node{(3,6,6)} \arrow{se,l,-}{2} \node[2]{(2,2,3,3)} \arrow{sw,r,-}{2} \\
  \node[2]{(3,3,3,3)}
  \end{diagram}
$$

\paragraph{\bf Class III}
$$
  \begin{diagram}
  \node{(2,6,8)} \arrow{s,l,-}{2} \arrow{ese,r,-}{2}
  \node[4]{(2,3,8)} \arrow{wsw,l,-}{10} \arrow{s,l,-}{2}
    \arrow{ese,r,-}{3^\ast} \\
  \node{(4,6,6)} \arrow{s,l,-}{3^\ast} \arrow{ese,l,-}{2}
  \node[2]{(3,8,8)} \arrow{s,l,-}{2}
  \node[2]{(3,3,4)} \arrow{wsw,l,-}{10} \arrow{s,r,-}{3^\ast}
  \node[2]{(2,4,8)} \arrow{wsw,r,-}{2} \arrow{s,r,-}{2} \\
  \node{(2^2,4^3)} \arrow{ese,l,-}{2}
  \node[2]{(3,4,3,4)} \arrow{s,l,-}{3^\ast}
  \node[2]{(4,4,4)} \arrow{wsw,l,-}{10} \arrow{s,r,-}{2}
  \node[2]{(2,8,8)} \arrow{wsw,l,-}{2} \arrow{s,r,-}{2} \\
  \node[3]{(4^6)} \arrow{s,l,-}{2}
  \node[2]{(2,4,2,4)} \arrow{wsw,l,-}{10} \arrow{s,r,-}{2}
  \node[2]{(4,8,8)} \arrow{wsw,r,-}{2} \\
  \node[3]{?} \arrow{s,l,-}{2}
  \node[2]{(4,4,4,4)} \arrow{wsw,r,-}{10} \\
  \node[3]{?}
  \end{diagram}
$$

\paragraph{\bf Class IV}
% $K=\Q(\sqrt 3)$ with $D=\fp_2$,
% $\Gamma^{(1)}=(3,3,6)$, $\Gamma^{(+)}=\Gamma^{(*)}=(2,3,12)$.
$$
  \begin{diagram}
  \node[3]{(2,3,12)} \arrow{wsw,l,-}{2} \arrow{s,l,-}{4} \arrow{ese,r,-}{3} \\
  \node{(3,3,6)} \arrow{s,l,-}{4} \arrow{ese,r,-}{} \node[2]{(3,4,12)}
    \arrow{wsw,r,-}{} \arrow{ese,r,-}{} \node[2]{(2,6,12)}
    \arrow{wsw,r,-}{} \arrow{s,r,-}{4} \arrow{ese,r,-}{2} \\
  \node{(3,3,2,6)} \arrow{ese,l,-}{3^\ast} \node[2]{(6,6,6)}
    \arrow{s,l,-}{4} \arrow{ese,r,-}{2} \node[2]{(2,4,6,12)}
    \arrow{wsw,l,-}{2} \arrow{ese,r,-}{} \node[2]{(3,12,12)} \arrow{wsw,r,-}{2}
    \arrow{s,r,-}{4} \\
  \node[3]{(2^3,6^3)} \arrow{ese,l,-}{2} \node[2]{(3^2,6^2)}
    \arrow{s,l,-}{4} \node[2]{(3,4^2,12^2)} \arrow{wsw,r,-}{2} \\
  \node[5]{(1;2,2,3,3,6,6)}
  \end{diagram}
$$

\paragraph{\bf Class V}
% $K=\Q(\sqrt 3)$ with $D=\fp_3$,
%$\Gamma^{(1)}=(2,2,2,6)$, $\Gamma^{(+)}=\Gamma^{(*)}=(2,4,12)$.
$$
  \begin{diagram}
  \node[3]{(2,4,12)} \arrow{sw,l,-}{2} \arrow{se,r,-}{2} \\
  \node[2]{(2,12,12)} \arrow{sw,l,-}{2} \arrow{se,r,-}{2}
  \node[2]{(4,4,6)} \arrow{sw,r,-}{2} \\
  \node{(6,12,12)} \arrow{se,l,-}{2} \node[2]{(2,2,6,6)} \arrow{sw,r,-}{2} \\
  \node[2]{(6,6,6,6)}
  \end{diagram}
$$

\paragraph{\bf Class VI}
% $K=\Q(\sqrt 5)$ with $D=\fp_2$,
%$\Gamma^{(1)}=\Gamma^{(+)}=(2,5,5)$, and $\Gamma^{(*)}=(2,4,5)$.
$$
  \begin{diagram}
  \node{(2,4,5)} \arrow{s,l,-}{2} \arrow{ese,r,-}{6}
    \node[4]{(2,4,10)} \arrow{wsw,l,-}{2} \arrow{s,r,-}{2} \\
  \node{(2,5,5)} \arrow{ese,l,-}{6^\ast}
    \node[2]{(4,4,5)} \arrow{s,l,-}{2}
    \node[2]{(2,10,10)} \arrow{wsw,l,-}{2} \arrow{s,r,-}{2} \\
  \node[3]{(2,2,5,5)} \arrow{s,l,-}{2}
    \node[2]{(5,10,10)} \arrow{wsw,r,-}{2} \\
  \node[3]{(5,5,5,5)}
  \end{diagram}
$$

\paragraph{\bf Class VII}
% $K=\Q(\sqrt 5)$ with $D=\fp_3$,
%$\Gamma^{(1)}=\Gamma^{(+)}=(3,5,5)$, and $\Gamma^{(*)}=(2,5,6)$.
$$
  \begin{diagram}
  \node{(2,5,6)} \arrow{s,r,-}{2} \\
  \node{(3,5,5)}
  \end{diagram}
$$

\paragraph{\bf Class VIII}
% $K=\Q(\sqrt 5)$ with $D=\fp_5$,
%$\Gamma^{(1)}=\Gamma^{(+)}=(3,3,5)$, and $\Gamma^{(*)}=(2,3,10)$.
$$
  \begin{diagram}
  \node[2]{(2,3,10)} \arrow{sw,l,-}{2} \arrow{se,r,-}{3} \\
  \node{(3,3,5)} \arrow{se,l,-}{3} \node[2]{(2,5,10)}
  \arrow{sw,r,-}{2} \\
  \node[2]{(5,5,5)}
  \end{diagram}
$$

\paragraph{\bf Class IX}
% $K=\Q(\sqrt 6)$ with $D=\fp_2$,
%$\Gamma^{(1)}=(2,3,3,3)$, and $\Gamma^{(+)}=\Gamma^{(*)}=(3,4,6)$.
$$
  \begin{diagram}
  \node{(3,4,6)}
  \end{diagram}
$$

\paragraph{\bf Class X}
% $K=\Q(\cos\pi/7)$ with $D=1$ and
%$\Gamma^{(1)}=\Gamma^{(+)}=\Gamma^{(*)}=(2,3,7)$.
$$
  \begin{diagram}
  \node{(2,4,7)} \arrow{se,l,-}{2}
    \node[2]{(2,3,7)} \arrow{sw,l,-}{9} \arrow{se,r,-}{8}
    \node[2]{(2,3,14)} \arrow{sw,r,-}{2} \arrow{se,r,-}{3} \\
  \node[2]{(2,7,7)} \arrow{se,l,-}{8}
    \node[2]{(3,3,7)} \arrow{sw,r,-}{9} \arrow{se,r,-}{3}
    \node[2]{(2,7,14)} \arrow{sw,r,-}{2} \\
  \node[3]{(1;7,7)} \arrow{se,r,-}{3}
  \node[2]{(7,7,7)} \arrow{sw,r,-}{9} \\
  \node[4]{(1;7^6)}
  \end{diagram}
$$

\paragraph{\bf Class XI}
% $K=\Q(\cos\pi/9)$ with $D=1$ and
%$\Gamma^{(1)}=\Gamma^{(+)}=\Gamma^{(*)}=(2,3,9)$.
%  \begin{diagram}
%  \node{(2,3,9)} \arrow{s,l,-}{4} \node[2]{(2,3,18)}
%  \arrow{wsw,l,-}{2} \arrow{s,r,-}{3} \arrow{ese,r,-}{4} \\
%  \node{(3,3,9)} \arrow{se,l,-}{3} \node[2]{(2,9,18)}
%  \arrow{sw,r,-}{2} \arrow{se,l,-}{4} \node[2]{(3,6,18)} \arrow{sw,r,-}{3} \\
%  \node[2]{(9,9,9)} \arrow{se,l,-}{4} \node[2]{(3,6,9,18)}
%  \arrow{sw,r,-}{2} \\
%  \node[3]{(3,3,3,9,9,9)}
%  \end{diagram}
%$$
%$$
%  \begin{diagram}
%  \node[3]{(2,3,18)} \arrow{wsw,l,-}{4} \arrow{s,l,-}{2} \arrow{ese,r,-}{3} \\
%  \node{(3,6,18)} \arrow{s,l,-}{2} \arrow{ese,r,-}{} \node[2]{(3,3,9)}
%    \arrow{wsw,r,-}{} \arrow{ese,r,-}{} \node[2]{(2,9,18)}
%    \arrow{wsw,r,-}{} \arrow{s,r,-}{2} \\
%  \node{(3,3,3,9)} \arrow{ese,l,-}{3^\ast} \node[2]{(3,6,9,18)}
%    \arrow{s,l,-}{2} \node[2]{(9,9,9)} \arrow{wsw,l,-}{4} \\
%  \node[3]{(3^3,9^3)}
%  \end{diagram}
%$$
$$
  \begin{diagram}
  \node{(2,3,9)} \arrow{s,l,-}{4}
  \node[2]{(2,3,18)} \arrow{wsw,l,-}{2} \arrow{s,l,-}{4} \arrow{ese,r,-}{3} \\
  \node{(3,3,9)} \arrow{s,l,-}{4} \arrow{ese,r,-}{} \node[2]{(3,6,18)}
    \arrow{wsw,r,-}{} \arrow{ese,r,-}{} \node[2]{(2,9,18)}
    \arrow{wsw,r,-}{} \arrow{s,r,-}{4} \\
  \node{(3,3,3,9)} \arrow{ese,l,-}{3^\ast} \node[2]{(9,9,9)}
    \arrow{s,l,-}{4} \node[2]{(3,6,9,18)} \arrow{wsw,l,-}{2} \\
  \node[3]{(3^3,9^3)}
  \end{diagram}
$$

\paragraph{\bf Class XII}
% $K=\Q(\cos\pi/9)$ with $D=\fp_2\fp_3$,
%$\Gamma^{(1)}=\Gamma^{(+)}=(0;2,2,9,9)$ and $\Gamma^{(\ast)}=(2,4,18)$.
$$
  \begin{diagram}
  \node[3]{(2,4,18)} \arrow{sw,l,-}{2} \arrow{se,r,-}{2} \\
  \node[2]{(2,18,18)} \arrow{sw,l,-}{2} \arrow{se,r,-}{2}
  \node[2]{(4,4,9)} \arrow{sw,r,-}{2} \\
  \node{(9,18,18)} \arrow{se,l,-}{2} \node[2]{(2,2,9,9)} \arrow{sw,r,-}{2} \\
  \node[2]{(9,9,9,9)}
  \end{diagram}
$$

\paragraph{\bf Class XIII}
% $K=\Q(\cos\pi/8)$ and $D=\fp_2$,
%$\Gamma^{(1)}=\Gamma^{(+)}=(3,3,8)$ and $\Gamma^{(\ast)}=(2,3,16)$.
$$
  \begin{diagram}
  \node[3]{(2,3,16)} \arrow{sw,l,-}{3} \arrow{se,r,-}{2} \\
  \node[2]{(2,8,16)} \arrow{sw,l,-}{2} \arrow{se,r,-}{2}
  \node[2]{(3,3,8)} \arrow{sw,r,-}{3} \\
  \node{(4,16,16)} \arrow{se,l,-}{2} \node[2]{(8,8,8)} \arrow{sw,r,-}{2} \\
  \node[2]{(4,4,8,8)}
  \end{diagram}
$$

\paragraph{\bf Class XIV}
% $K=\Q(\cos\pi/10)$ and $D=\fp_2$,
%$\Gamma^{(1)}=(5,5,10)$, $\Gamma^{(+)}=\Gamma^{(\ast)}=(2,5,20)$.
$$
  \begin{diagram}
  \node{(2,5,20)} \arrow{s,r,-}{2} \\
  \node{(5,5,10)}
  \end{diagram}
$$

\paragraph{\bf Class XV}
% $K=\Q(\cos\pi/12)$ and $D=\fp_2$,
%$\Gamma^{(1)}=(3,3,12)$, $\Gamma^{(+)}=\Gamma^{(\ast)}=(2,3,24)$.
%$$
%  \begin{diagram}
%  \node[3]{(2,3,24)} \arrow{wsw,l,-}{4} \arrow{s,l,-}{2} \arrow{ese,r,-}{3} \\
%  \node{(3,8,24)} \arrow{se,l,-}{2} \node[2]{(3,3,12)} \arrow{sw,l,-}{4}
%    \arrow{se,r,-}{3} \node[2]{(2,12,24)} \arrow{sw,r,-}{2}
%    \arrow{se,r,-}{2} \\
%  \node[2]{(3,3,4,12)} \arrow{se,l,-}{3^\ast} \node[2]{(12,12,12)}
%    \arrow{sw,l,-}{4} \arrow{se,r,-}{2} \node[2]{(6,24,24)}
%    \arrow{sw,r,-}{2} \\
%  \node[3]{(4^3,12^3)} \arrow{se,l,-}{} \node[2]{(6^2,12^2)}
%    \arrow{sw,r,-}{} \\
%  \node[4]{?}
%  \end{diagram}
%$$
$$
  \begin{diagram}
  \node[3]{(2,3,24)} \arrow{wsw,l,-}{2} \arrow{s,l,-}{4} \arrow{ese,r,-}{3} \\
  \node{(3,3,12)} \arrow{s,l,-}{4} \arrow{ese,r,-}{} \node[2]{(3,8,24)}
    \arrow{wsw,r,-}{} \arrow{ese,r,-}{} \node[2]{(2,12,24)}
    \arrow{wsw,r,-}{} \arrow{s,r,-}{4} \arrow{ese,r,-}{2} \\
  \node{(3,3,4,12)} \arrow{ese,l,-}{3^\ast} \node[2]{(12,12,12)}
    \arrow{s,l,-}{4} \arrow{ese,r,-}{2} \node[2]{(4,8,12,24)}
    \arrow{wsw,l,-}{2} \arrow{ese,r,-}{} \node[2]{(6,24,24)} \arrow{wsw,r,-}{2}
    \arrow{s,r,-}{4} \\
  \node[3]{(4^3,12^3)} \arrow{ese,l,-}{2} \node[2]{(6^2,12^2)}
    \arrow{s,l,-}{4} \node[2]{(2,6,8^2,24^2)} \arrow{wsw,r,-}{2} \\
  \node[5]{(1;2^2,4^2,6^2,12^2)}
  \end{diagram}
$$

\paragraph{\bf Class XVI}
% $K=\Q(\cos\pi/15)$ and $D=\fp_3$,
%$\Gamma^{(1)}=(5,5,15)$, $\Gamma^{(+)}=\Gamma^{(\ast)}=(2,5,30)$.
$$
  \begin{diagram}
  \node{(2,5,30)} \arrow{s,r,-}{2} \\
  \node{(5,5,15)}
  \end{diagram}
$$

\paragraph{\bf Class XVII}
% $K=\Q(\cos\pi/15)$ and $D=\fp_5$,
%$\Gamma^{(1)}=(3,3,15)$, $\Gamma^{(+)}=\Gamma^{(\ast)}=(2,3,30)$.
$$
  \begin{diagram}
  \node[3]{(2,3,30)} \arrow{wsw,l,-}{2} \arrow{s,l,-}{4} \arrow{ese,r,-}{3} \\
  \node{(3,3,15)} \arrow{s,l,-}{4} \arrow{ese,r,-}{} \node[2]{(3,10,30)}
    \arrow{wsw,r,-}{} \arrow{ese,r,-}{} \node[2]{(2,15,30)}
    \arrow{wsw,r,-}{} \arrow{s,r,-}{4} \\
  \node{(3,3,5,15)} \arrow{ese,l,-}{3^\ast} \node[2]{(15,15,15)}
    \arrow{s,l,-}{4} \node[2]{(5,10,15,30)} \arrow{wsw,l,-}{2} \\
  \node[3]{(5^3,15^3)}
  \end{diagram}
$$

\paragraph{\bf Class XVIII}
% $K=\Q(\sqrt 2,\sqrt 5)$ and $D=\fp_2$,
%$\Gamma^{(1)}=\Gamma^{(+)}=(4,5,5)$, and $\Gamma^{(\ast)}=(2,5,8)$.
$$
  \begin{diagram}
  \node{(2,5,8)} \arrow{s,r,-}{2} \\
  \node{(4,5,5)}
  \end{diagram}
$$

\paragraph{\bf Class XIX}
% $K=\Q(\cos\pi/11)$ with $D=1$ and
%$\Gamma^{(1)}=\Gamma^{(+)}=\Gamma^{(\ast)}=(2,3,11)$.
$$
  \begin{diagram}
  \node{(2,3,11)}
  \end{diagram}
$$

\end{document}